\documentclass[11pt]{article}
\usepackage[tbtags]{amsmath}
\usepackage{amssymb}
\usepackage{amsthm}
\usepackage[misc]{ifsym}
\usepackage{cases}
\usepackage{mathrsfs}
\numberwithin{equation}{section}   %????????equation?????§ß?????section ????
\normalsize
\title{\bf {An overlapping information linear-quadratic Stackelberg stochastic differential game with two leaders and two followers} \thanks{This work is supported by National Key R\&D Program of China (Grant No. 2022YFA1006104), National Natural Science Foundations of China (Grant Nos. 11971266, 12271304, 11831010), and Shandong Provincial Natural Science Foundations (Grant Nos. ZR2022JQ01, ZR2020ZD24, ZR2019ZD42).}}
\author{\normalsize  Yu Si\thanks{\it School of Mathematics, Shandong University, Jinan 250100, P.R. China, E-mail: 202112003@mail.sdu.edu.cn} , Jingtao Shi\thanks{\it Corresponding author. School of Mathematics, Shandong University, Jinan 250100, P.R. China, E-mail: shijingtao@sdu.edu.cn}}
\newtheorem{Remark}{Remark}[section]
\newtheorem{assumption}{Assumption}[section]
\begin{document}
\maketitle

\noindent{\bf Abstract:}\quad This paper is concerned with an overlapping information linear-quadratic (LQ) Stackelberg stochastic differential game with two leaders and two followers, where the diffusion terms of the state equation contain both the control and state variables. A distinct feature lies in that, the noisy information available to the leaders and the followers may be asymmetric and have overlapping part. Using a coupled system of Riccati equations, the followers first solve an LQ nonzero-sum stochastic differential Nash game with partial information, and then the leaders solve a partial information LQ nonzero-sum stochastic differential Nash game driven by a conditional mean-field type forward-backward stochastic differential equation (CMF-FBSDE). By maximum principle, completion of squares and decoupling methods, the state-estimate feedback representation of the Stackelberg-Nash equilibrium is obtained.

\vspace{2mm}


\vspace{2mm}


\section{Introduction}

In recent years, the Stackelberg (also known as leader-follower) game has become an active topic in noncooperative game research. Its characteristic is that the decisions of the two players are made in sequence, with the leader first declaring his action, and the follower optimizing his cost functional in response. The leader then seeks her strategy based on the follower's rational response to optimize her cost functional. Most of the research on Stackelberg games only involves the case of a single leader. While in real life, there are many situations with multiple leaders, such as the collective leadership system in political science, whose main feature is that decision-making is collectively or institutionally responsible, rather than being made by one person. Another example exists in the military, the commander is responsible for military affairs, while the political commissar is responsible for political affairs, which presents a game problem with two leaders. Based on these real-life examples, this paper proposes a Stackelberg game problem with two leaders and two followers with overlapping information.

The research of Stackelberg games can be traced back to the pioneering work by Stackelberg \cite{Stackelberg-1952} in static competitive economics. There is a lengthy body of literature on Stackelberg games, which we will not elaborate on. Let us only mention some of them related to the topic of this paper. Simann and Cruz \cite{Simaan-Cruz-1976} studied the dynamic LQ Stackelberg differential game with two players, and the Stackelberg strategy was expressed in terms of Riccati-like differential equations. Yong \cite{Yong-2002} researched the LQ Stackelberg stochastic differential game to a rather general framework, where the coefficients could be random matrices, the control variables could enter the diffusion term of the state equation and the weight matrices for the controls in the cost functionals need not to be positive definite. Open-loop Stackelberg equilibrium and its state feedback representation of two players is obtained, firstly, by maximum principle and theory of FBSDE. \O ksendal et al. \cite{Oksendal-Sandal-Uboe-2013} proved a maximum principle for the Stackelberg differential game with jumps, and applied the result to a continuous-time manufacturer-newsvendor model. Bensoussan et al. \cite{Bensoussan-Chen-Sethi-2015} proposed several solution concepts about information structure for the stochastic Stackelberg differential game with two players, and derived the maximum principle under closed-loop memoryless information structure. Lin et al. \cite{Lin-Jiang-Zhang-2019} investigated an open-loop Stackelberg strategy for the LQ stochastic differential game of two players, where the game system is governed by a mean-field stochastic differential equation (MF-SDE). Very recently, Sun et al. \cite{Sun-Wang-Wen-2023} studied a two-player zero-sum LQ Stackelberg stochastic differential game, and the Stackelberg equilibrium was obtained by first solving a stochastic LQ optimal control problem and then a backward stochastic LQ optimal control problem.

There are some literatures about Stackelberg stochastic differential games with multiple followers. Mukaidani and Xu \cite{Mukaidani-Xu-2015} studied a Stackelberg stochastic differential game with one leader and multiple followers. Stackelberg strategies were developed in the settings that the followers act either cooperatively to attain Pareto optimality or non-cooperatively to arrive at a Nash equilibrium. Li and Yu \cite{Li-Yu-2018} used a kind of FBSDEs with a self-similar domination-monotonicity structure to characterize the unique equilibrium of an LQ generalized Stackelberg stochastic differential game with multilevel hierarchy. Wang and Zhang \cite{Wang-Zhang-2020} studied an LQ Stackelberg stochastic differential game of mean-field type, with one leader and two followers. By maximum principle and verification theorem, the open-loop Stackelberg solution was expressed as a feedback form of the state and its mean with the help of three systems of Riccati equations. Wang and Yan \cite{Wang-Yan-2023} obtained the Pareto-based Stackelberg equilibriuma for a Stackelberg stochastic differential game with multiple followers.

In recent years, more and more research attention is being drawn to the theory of stochastic differential games with asymmetric and overlapping information, which has wide applicable background in, for example, insider tradings (Biagini and \O ksendal \cite{Biagini-Oksendal-2006}), principal-agent/optimal contract problems (Cvitanic and Zhang \cite{Cvitanic-Zhang-2013}). Chang and Xiao \cite{Chang-Xiao-2014} studied an LQ nonzero-sum differential game with asymmetric information. Nash equilibria are obtained for several classes of asymmetric information by stochastic maximum principle (SMP) and represented in a feedback form of the optimal filtering of the state, through the solutions of the Riccati equations. Shi et al. \cite{Shi-Wang-Xiong-2016} introduced a general framework of Stackelberg stochastic differential game with asymmetric information, and derived the SMP and verification theorem with partial information, to represent the Stackelberg equilibrium. Shi et al. \cite{Shi-Wang-Xiong-2017,Shi-Wang-Xiong-2020} researched LQ Stackelberg stochastic differential games with asymmetric information and overlapping information, respectively. Wang et al. \cite{Wang-Wang-Zhang-2020} investigated a Stackelberg stochastic differential game with asymmetric information of one-leader and two-followers. The open-loop Stackelberg solution was expressed as a feedback form of state, its estimation and mean. Zhao et al. \cite{Zhao-Zhang-Gao-Lv-2022} researched an incomplete information LQ Stackelberg stochastic differential game of two leaders and two followers, under the hierarchical network architecture of industrial internet of things (IIoT). The feedback Stackelberg-Nash equilibrium point and the corresponding dynamic equations of the game were derived by using the SMP. Zheng and Shi \cite{Zheng-Shi-2022} introduced a general setting of Stackelberg stochastic differential games with asymmetric noisy observation, and LQ case was studied via the SMP with partial information and that with conditional mean-field type FBSDE (CMF-FBSDE). Very recently, Kang and Shi \cite{Kang-Shi-2023} studied a three-level LQ Stackelberg stochastic differential game with asymmetric information. By the SMP of FBSDEs and optimal filtering, feedback Stackelberg equilibrium was obtained with a system of three Riccati equations.

Besides what are mentioned above, in some cases we may encounter multiple leaders in Stackelberg games in reality. However, as far as we know, there a significant lack of literature in this research field. Sherali \cite{Sherali-1984} introduced a multiple leader-follower model which extended Stackelberg \cite{Stackelberg-1952}'s classical one to the duopoly. DeMiguel and Xu \cite{DeMiguel-Xu-2009} studied an oligopolistic Stacklelberg stochastic game consisting of multiple leader-follower to supply a homogeneous product (or service) noncooperatively, and obtained the Stackelberg-Nash-Cournot (SNC) equalibrium. Huang et al. \cite{Huang-Si-Wu-2021} studied an LQ stochastic large population system combining three types of interactive agents mixed, which are respectively, major leader, minor leaders, and minor followers. The SNC approximate equilibrium was derived from the combination of a major-minor mean-field game (MFG) and a leader-follower Stackelberg game. Zhao et al. \cite{Zhao-Zhang-Gao-Lv-2022} researched an incomplete information LQ Stackelberg stochastic differential game of multiple leader-follower, under the hierarchical network architecture of IIoT.

Motivated by the above literature and practical applications, in this paper, we consider an overlapping information LQ Stackelberg differential game with multiple leader-follower. We only consider the model of two leaders and two followers. Compared to \cite{Zhao-Zhang-Gao-Lv-2022}, the information of the two followers is the same, and the information of the two leaders is identical. However, there is overlapping information between followers and leaders. Additionally, the state and control variables enter the system's diffusion terms. These new features make the problem more difficult and challenging. We divide the game problem into two parts. The followers meet an LQ two-player nonzero-sum Nash game driven by SDE, and the leaders need to solve an LQ two-player nonzero-sum Nash game driven by a CMF-FBSDE. To overcome the difficulty caused by the overlapping information, we use the SMP with partial information to get the open-loop form of the Stackelberg-Nash equilibrium, and then using the methods of undetermined coefficients and dimension expansion to obtain its state-estimate feedback representation.

The rest of this paper is organized as follows. In Section 2, we formulate our problem. In Section 3, we address the problem of the followers and leaders in turn and derive the main results. Finally, conclusion is given in Section 4.

\section{Problem formulation}

Let $\left(\Omega, \mathcal{F}, \mathbb{P},\left\{\mathcal{F}_t\}_{t\ge0}\right)\right.$ be a complete filtered probability space, on which a standard three-dimensional Brownian motion $\left(W_1(\cdot), W_2(\cdot), W_3(\cdot)\right)$ is defined, $\left\{\mathcal{F}_t\right\}_{0 \leq t \leq T}$ is the natural filtration generated by $\left(W_1(\cdot), W_2(\cdot), W_3(\cdot)\right)$, and $T>0$ is a fixed time duration. Suppose that the state process $x^{u_1, u_2, u_3, u_4}(\cdot)$ satisfies a linear SDE:
\begin{equation}\label{state}
\left\{\begin{aligned}
d x^{u_1, u_2, u_3, u_4}(t) & =\Big[a_0(t) x^{u_1, u_2, u_3, u_4}(t)+b_0(t) u_1(t)+c_0(t) u_2(t)+d_0(t) u_3(t) \\
&\quad +e_0(t) u_4(t)\Big] d t + \sum_{i=1}^3 \Big[a_i(t) x^{u_1, u_2, u_3, u_4}(t)+b_i(t) u_1(t)+c_i(t) u_2(t) \\
&\quad +d_i(t) u_3(t)+e_i(t) u_4(t)\Big] d W_i(t),\quad t\in[0,T], \\
x^{u_1, u_2, u_3, u_4}(0) & =x_0,
\end{aligned}\right.
\end{equation}
where $a_i(\cdot), b_i(\cdot), c_i(\cdot), d_i(\cdot), e_i(\cdot)$ are $\mathbb{R}$-valued deterministic and bounded functions, $i=0,1,2,3$. Here, $x^{u_1, u_2, u_3, u_4}(\cdot) \in \mathbb{R}$ is the state process, and $u_j(\cdot) \in \mathbb{R}$ is the control process of player $j$, for $j=1,2,3,4$. In our game of this paper, we suppose that players 1, 2 are the followers, and players 3, 4 are the leaders. Moreover, there is asymmetric information between the followers and leaders. We define $\mathcal{G}_t^1\equiv\mathcal{F}_t^j=\sigma\left\{W_1(s), W_3(s) ; 0 \leq s \leq t\right\}, j=1,2$ which denotes the information available to the followers and $\mathcal{G}_t^2\equiv\mathcal{F}_t^j=\sigma\left\{W_2(s), W_3(s) ; 0 \leq s \leq t\right\}, j=3,4$ to the leaders. Then we can define the admissible control sets of the followers and the leaders:
$$
\mathcal{U}_j:=\Big\{u_j(\cdot) \mid u_j(\cdot): \Omega \times[0, T] \rightarrow \mathbb{R} \text { is } \mathcal{F}_t^j \text {-adapted and square integrable }\Big\},\quad j=1,2,3,4.
$$
For any $u_j(\cdot) \in \mathcal{U}_j $, $j=1,2,3,4$, it is classical that SDE (\ref{state}) admits a unique $\mathcal{F}_t$-adapted solution $x^{u_1, u_2, u_3, u_4}(\cdot)$. For $j=1,2,3,4$, we define the cost functional of the player $j$ as
\begin{equation}
\begin{aligned}\label{cost}
&J_j(u_1(\cdot), u_2(\cdot), u_3(\cdot), u_4(\cdot))\\
&= \frac{1}{2} \mathbb{E}\left\{\int_0^T\Big[ Q_j(t) x^{u_1, u_2, u_3, u_4}(t)^2+ R_j(t) u_j(t)^2\Big] d t+ G_j x^{u_1, u_2, u_3, u_4}(T)^2\right\},
\end{aligned}
\end{equation}
where $Q_j(\cdot), R_j(\cdot)$ are $\mathbb{R}$-valued deterministic and bounded functions, and $G_j$ are constants.

Now, let us formulate the Stackelberg game among the two leaders and two followers. First, for given $\left(u_3(\cdot), u_4(\cdot)\right)\in \mathcal{U}_3 \times \mathcal{U}_4$ of the leaders, the two followers face a nonzero-sum Nash game, that is, they choose $\left(u_1^*(\cdot), u_2^*(\cdot)\right) \equiv\left(u_1^*\left(\cdot, u_3(\cdot), u_4(\cdot)\right), u_2^*\left(\cdot, u_3(\cdot), u_4(\cdot)\right)\right) \in$ $\mathcal{U}_1 \times \mathcal{U}_2$, such that
\begin{equation}\label{Nash game of followers}
\left\{\begin{array}{l}
J_1\left(u_1^*(\cdot), u_2^*(\cdot), u_3(\cdot), u_4(\cdot)\right) = \inf\limits_{u_1(\cdot) \in\, \mathcal{U}_1} J_1\left(u_1(\cdot), u_2^*(\cdot), u_3(\cdot),  u_4(\cdot)\right), \\
J_2\left(u_1^*(\cdot), u_2^*(\cdot), u_3(\cdot), u_4(\cdot)\right) =\inf\limits_{u_2(\cdot) \in\, \mathcal{U}_2} J_2\left(u_1^*(\cdot), u_2(\cdot), u_3(\cdot), u_4(\cdot)\right) .
\end{array}\right.
\end{equation}
Then, knowing that the followers would take $\left(u_1^*(\cdot), u_2^*(\cdot)\right)$, the two leaders are also faced with a nonzero-sum Nash game, that is, they would like to find $\left(u_3^*(\cdot), u_4^*(\cdot)\right)$, such that
\begin{equation}\label{Nash game of leaders}
\left\{\begin{array}{l}
J_3\left(u_1^*(\cdot), u_2^*(\cdot), u_3^*(\cdot), u_4^*(\cdot)\right) = \inf\limits_{u_3(\cdot) \in\, \mathcal{U}_3} J_3\left(u_1^*(\cdot), u_2^*(\cdot), u_3(\cdot), u_4^*(\cdot)\right), \\
J_4\left(u_1^*(\cdot), u_2^*(\cdot), u_3^*(\cdot), u_4^*(\cdot)\right) =\inf\limits_{u_4(\cdot) \in\, \mathcal{U}_4} J_4\left(u_1^*(\cdot), u_2^*(\cdot), u_3^*(\cdot), u_4(\cdot)\right) .
\end{array}\right.
\end{equation}
Noting that the information available to the leaders and followers has overlapping part, we refer to the above problem as {\it an overlapping information LQ Stackelberg stochastic differential game with two leaders and two followers}. If such an optimal control quadruple $(u_1^*(\cdot), u_2^*(\cdot), u_3^*(\cdot), u_4^*(\cdot))$ exists, we call it an {\it open-loop Stackelberg-Nash equilibrium} of the game.

\section{Main result}

In this section, we deal with the problems of the follower and the leader in two subsections, respectively.
For any $\mathcal{F}_t$-adapted process $\xi(\cdot)$, we denote by
$$
\hat{\xi}(t):=\mathbb{E}\left[\xi(t) | \mathcal{G}_t^1\right],\quad \check{\xi}(t):=\mathbb{E}\left[\xi(t) | \mathcal{G}_t^2\right],\quad
\check{\hat{\xi}}(t):=\mathbb{E}\left[\mathbb{E}\left[\xi(t) | \mathcal{G}_t^1\right] | \mathcal{G}_t^2\right] \equiv \mathbb{E}\left[\mathbb{E}\left[\xi(t) | \mathcal{G}_t^2\right] | \mathcal{G}_t^1\right]
$$
its optimal filtering estimates, for $t\in[0,T]$. Hereinafter we often drop the time variable $t$ for simplicity if no confusion arises.

\subsection{Problem of the followers}

For given $\left(u_3(\cdot), u_4(\cdot)\right)\in \mathcal{U}_3 \times \mathcal{U}_4$, the followers are facing a nonzero-sum stochastic differential Nash game with partial information.

First, we introduce the following assumption.
\begin{assumption}\label{A1}
	$Q_j(t) \geq 0, R_j(t) \geq 0 ,\forall t \in[0, T]$ and $G_j \geq 0$, for $j=1,2$.
\end{assumption}

We have the following result and its proof is postponed in the Appendix.
\newtheorem{thm}{Theorem}[section]
\begin{thm}\label{thm1}
Under Assumption \ref{A1}, $\left(u_1^*(\cdot), u_2^*(\cdot)\right)$ is a Nash equilibrium point of the followers' problem, if and only if
\begin{equation}\label{SMP of followers}
			R_1 u_1^*=b_0 \hat{p}_1+\sum\limits_{i=1}^3 b_i \hat{k}_{1 i}, \quad
			R_2 u_2^*=c_0 \hat{p}_2+\sum\limits_{i=1}^3 c_i \hat{k}_{2 i},
\end{equation}
where $\left(p_j(\cdot), k_{j1}(\cdot), k_{j2}(\cdot), k_{j3}(\cdot) \right)$, $j=1,2$ are the unique $\mathcal{F}_t$-adapted solutions satisfying BSDEs:
\begin{equation}\label{adjoint BSDEs of followers}
		\left\{\begin{aligned}
			-dp_1(t)&=\bigg[a_0 p_1+\sum\limits_{i=1}^3 a_i k_{1 i}-Q_1 x^{u_{1}^*, u_2^*, u_3, u_4}\bigg] d t-\sum\limits_{i=1}^3 k_{1 i} d W_i , \\
			-dp_2(t)&=\bigg[a_0 p_2+\sum\limits_{i=1}^3 a_i k_{2 i}-Q_2 x^{u_1^* u_2^*, u_3, u_4}\bigg] d t-\sum\limits_{i=1}^3 k_{2 i} d W_i, \quad t\in[0,T],\\
			  p_1(T)&=-G_1 x^{u_1^*, u_2^*, u_3, u_4}(T), \quad p_2(T)=-G_2 x^{u_1^*, u_2^*, u_3, u_4}(T).
		\end{aligned}
		\right.
\end{equation}
\end{thm}

Owing to that (\ref{adjoint BSDEs of followers}) is hard to apply, we will derive a state-estimate feedback form of the Nash equilibrium point $\left(u_1^*(\cdot), u_2^*(\cdot)\right)$.
For this target, we write the following optimality system:
\begin{equation}
\left\{\begin{aligned}\label{optimality system of followers}
d x^{u_3, u_4}(t) & =\Big[a_0 x^{u_3, u_4}+b_0 u_1^*+c_0 u_2^*+d_0 u_3+e_0 u_4\Big] d t  \\
&\quad  + \sum_{i=1}^3 \Big[a_i x^{u_3, u_4}+b_i u_1^*+c_i u_2^*+d_i u_3+e_i u_4\Big] d W_i,   \\
- d p_1(t) & =\bigg[a_0 p_1+\sum\limits_{i=1}^3 a_i k_{1 i}-Q_1 x^{u_3, u_4}\bigg] d t-\sum\limits_{i=1}^3 k_{1 i} d W_i , \\
- d p_2(t) & =\bigg[a_0 p_2+\sum\limits_{i=1}^3 a_i k_{2 i}-Q_2 x^{u_3, u_4}\bigg] d t-\sum_{i=1}^3 k_{2 i} d W_i, \quad t\in[0,T],\\
x^{u_3, u_4}(0)  & = x_0, \\
p_1(T) &= -G_1 x^{u_3, u_4}(T),\quad p_2(T) = -G_2 x^{u_3, u_4}(T),
\end{aligned}\right.
\end{equation}
where $x^{u_3, u_4}(\cdot)\equiv x^{u_1^*, u_2^*, u_3, u_4}(\cdot)$ for simplicity.

Observing the terminal condition of (\ref{optimality system of followers}), for $j=1,2$, we put
\begin{equation}\label{relation between adjoint and state of followers}
p_j(t)=-\tilde{P}_j(t) x^{u_3, u_4}(t)-\phi_j(t),\quad t\in[0,T],
\end{equation}
where $\tilde{P}_j(\cdot) \in \mathbb{R}$ is a deterministic and differentiable function, $\phi_j(\cdot) \in \mathbb{R}$ is an $\mathcal{F}_t$-adapted process and satisfies the BSDE
\begin{equation}\label{variphi BSDEs of followers}
\left\{\begin{array}{l}
d \phi_j(t)=\alpha_j(t) d t+\zeta_{j1}(t) d W_1+\zeta_{j3}(t) d W_3, \quad t\in[0,T],\\
\phi_j(T)=0,
\end{array}\right.
\end{equation}
where $\mathcal{F}_t$-adapted processes $\alpha_j(\cdot) \in \mathbb{R}$ and $\zeta_{j1}(\cdot), \zeta_{j3}(\cdot) \in \mathbb{R}$ will be determined later. Applying It\^{o}'s formula to $p_j(\cdot)$ in (\ref{relation between adjoint and state of followers}), we get
\begin{equation}\label{Ito formula of followers}
\begin{aligned}
d p_j &= -\dot{\tilde{P}}_j x^{u_3, u_4} d t-\tilde{P}_j\big(a_0 x^{u_3, u_4}+b_0 u_1^*+c_0 u_2^*+d_0 u_3+e_0 u_4\big) d t -\alpha_j d t\\
&\quad -\tilde{P}_j \sum_{i=1}^3\Big(a_i x^{u_3, u_4}+b_i u_1^*+c_i u_2^*+d_i u_3+e_i u_4\Big) d W_i  - \zeta_{j1} d W_1-\zeta_{j3} d W_3 \\
&= -\bigg(a_0 p_j+\sum_{i=1}^{3} a_i k_{j i}-Q_j x^{ u_3, u_4}\bigg) d t+\sum_{i=1}^3 k_{j i} d W_i, \quad j=1,2.
\end{aligned}
\end{equation}
Comparing the terms on both sides of the backward equations in (\ref{Ito formula of followers}), we obtain for $j=1,2$,
\begin{equation}\label{comparing of followers}
\left\{\begin{aligned}
-\alpha_j&= \dot{\tilde{P}}_j x^{u_3 , u_4}+\tilde{P}_j\big(a_0 x^{u_3, u_4}+b_0 u_1^*+c_0 u_2^*+d_0 u_3+e_0 u_4\big) -a_0 p_j-\sum_{i=1}^3 a_i k_{j i}+Q_j x^{u_3 u_4}, \\
k_{j i}&= -\tilde{P}_j\big(a_i x^{u_3, u_4}+b_i u_1^*+c_i u_2^*+d_i u_3+e_i u_4\big)-\zeta_{j i},\quad i=1,3, \\
k_{j2}&= -\tilde{P}_j\big(a_2 x^{u_3, u_4}+b_2 u_1^*+c_2 u_2^*+d_2 u_3+e_2 u_4\big).
\end{aligned}\right.
\end{equation}

Taking $\mathbb{E}\left[\cdot \mid \mathcal{G}_t^1\right]$ on both sides of (\ref{relation between adjoint and state of followers}) and (\ref{comparing of followers}), we have for $j=1,2$,
\begin{equation}\label{filtering relation of followers}
\begin{aligned}
\hat{p}_j & =-\tilde{P}_j \hat{x}^{u_3, u_4}-\hat{\phi}_j,
\end{aligned}
\end{equation}
and
\begin{equation}\label{filtering comparing of followers}
\left\{\begin{aligned}
-\hat{\alpha}_j &= \dot{\tilde{P}}_j \hat{x}^{u_3 , u_4}+\tilde{P}_j\big(a_0 \hat{x}^{u_3, u_4}+b_0 u_1^*+c_0 u_2^*+d_0 \hat{u}_3+e_0 \hat{u}_4\big)-a_0 \hat{p}_j-\sum_{i=1}^3 a_i \hat{k}_{j i}+Q_j \hat{x}^{u_3 u_4}, \\
\hat{k}_{j i} &= -\tilde{P}_j\big(a_i \hat{x}^{u_3, u_4}+b_i u_1^*+c_i u_2^*+d_i \hat{u}_3+e_i \hat{u}_4\big)-\hat{\zeta}_{j i},\quad i=1,3, \\
\hat{k}_{j2} &= -\tilde{P}_j\big(a_2 \hat{x}^{u_3, u_4}+b_2 u_1^*+c_2 u_2^*+d_2 \hat{u}_3+e_2 \hat{u}_4\big).
\end{aligned}\right.
\end{equation}
Substitute (\ref{filtering relation of followers}) and (\ref{filtering comparing of followers}) into (\ref{SMP of followers}), we achieve
\begin{equation}\label{SMP-SMP of followers}
\begin{aligned}
\left(R_1+\sum_{i=1}^3 \tilde{P}_1 b_i^2\right) u_1^*= & -\left\{\left(b_0 \tilde{P}_1+\sum_{i=1}^3 b_i \tilde{P}_1 a_i\right) \hat{x}^{u_3, u_4}\right.\left.\right.  \\
& \left.+\sum_{i=1}^3 b_i \tilde{P}_1\left(c_i u_2^*+d_i \hat{u}_3+e_i \hat{u}_4\right) +b_1 \hat{\zeta}_{11}+b_3 \hat{\zeta}_{13}+b_0 \hat{{\phi}}_1\right\}, \\
\left(R_2+\sum_{i=1}^3 \tilde{P}_2 c_i^2\right) u_2^*= & -\left\{\left(c_0 \tilde{P}_2+\sum_{i=1}^3 c_i \tilde{P}_2 a_i\right) \hat{x}^{u_3, u_4}\right.   \\
& \left.+\sum_{i=1}^3 c_i \tilde{P}_2\left(b_i u_1^*+d_i \hat{u}_3+e_i \hat{u}_4\right) +c_1 \hat{\zeta}_{21}+c_3 \hat{\zeta}_{23}+c_0 \hat{\phi}_2\right\}. \\
\end{aligned}
\end{equation}

We wish to solve $u_1^*$ and $u_2^*$ explicitly from the above coupled system (\ref{SMP-SMP of followers}). Set
\begin{equation*}
\begin{aligned}
 \bar{R}_1:= R_1+\sum_{i=1}^3 \tilde{P}_1 b_i^2,\quad \bar{R}_2:= R_2+\sum_{i=1}^3 \tilde{P}_2 c_i^2,\quad \bar{B}_1:= b_0+\sum_{i=1}^3 a_i b_i,\quad \bar{C}_1:= c_0+\sum_{i=1}^3 a_i c_i,
\end{aligned}
\end{equation*}
and we can obtain
{\small$$
\begin{aligned}
& \left(\begin{array}{cc}
\bar{R}_1 & \tilde{P}_1 \sum\limits_{i=1}^3 b_i c_i \\
\tilde{P}_2 \sum\limits_{i=1}^3 b_i c_i & \bar{R}_2
\end{array}\right)\left(\begin{array}{c}
u_1^* \\
u_2^*
\end{array}\right)=-\left[\left(\begin{array}{c}
\tilde{P}_1\bar{B}_1 \\
\tilde{P}_2\bar{C}_1
\end{array}\right) \hat{x}^{u_3, u_4}+\left(\begin{array}{c}
\tilde{P}_1\sum\limits_{i=1}^3 b_i d_i \\
\tilde{P}_2\sum\limits_{i=1}^3 c_i d_i
\end{array}\right) \hat{u}_3\right. \\
& \left.+\left(\begin{array}{c}
\tilde{P}_1\sum\limits_{i=1}^3 b_i e_i \\
\tilde{P}_2\sum\limits_{i=1}^3 c_i e_i
\end{array}\right) \hat{u}_4+\left(\begin{array}{cc}
b_0 & 0 \\
0 & c_0
\end{array}\right)\left(\begin{array}{c}
\hat{\phi}_1 \\
\hat{\phi}_2
\end{array}\right)+\left(\begin{array}{cc}
b_1 & 0 \\
0 & c_1
\end{array}\right)\left(\begin{array}{c}
\hat{\zeta}_{11} \\
\hat{\zeta}_{21}
\end{array}\right)+\left(\begin{array}{cc}
b_3 & 0 \\
0 & c_3
\end{array}\right)\left(\begin{array}{c}
\hat{\zeta}_{13} \\
\hat{\zeta}_{23}
\end{array}\right)\right].
\end{aligned}
$$}
Put
\begin{equation*}
\begin{aligned}
& N:=\left(\begin{array}{ll}
R_1+\sum\limits_{i=1}^3 \tilde{P}_1 b_i^2 & \tilde{P}_1 \sum\limits_{i=1}^3 b_i c_i \\
\tilde{P}_2\sum\limits_{i=1}^3 b_i c_i & R_2+\sum\limits_{i=1}^3 \tilde{P}_2 c_i^2
\end{array}\right)\equiv\left(\begin{array}{cc}
\bar{R}_1 & \tilde{P}_1 \sum\limits_{i=1}^3 b_i c_i \\
\tilde{P}_2 \sum\limits_{i=1}^3 b_i c_i & \bar{R}_2
\end{array}\right).
\end{aligned}
\end{equation*}
We need to introduce the following condition.
\begin{assumption}\label{A2}
	$N(t)$ is invertible for all $t \in [0, T]$.
\end{assumption}
%Then
%$$
%\begin{aligned}
%& N^{-1}=\frac{1}{|N|}\left(\begin{array}{cc}
%\bar{R}_2 & -\tilde{P}_1 \sum\limits_{i=1}^3 b_i c_i \\
%-\bar{P}_2 \sum\limits_{i=1}^3 b_i c_i & \bar{R}_1
%\end{array}\right).
%\end{aligned}
%$$
Therefore, by (\ref{SMP-SMP of followers}),
$$
\begin{aligned}
&\begin{aligned}
\left(\begin{array}{l}
u_1^* \\
u_2^*
\end{array}\right)=& -\frac{1}{|N|}\left(\begin{array}{cc}
\bar{R}_2 & -\tilde{P}_1 \sum\limits_{i=1}^3 b_i c_i \\
-\tilde{P}_2 \sum\limits_{i=1}^3 b_i c_i & \bar{R}_1
\end{array}\right)\Bigg[\left(\begin{array}{c}
\tilde{P}_1\bar{B}_1 \\
\tilde{P}_2\bar{C}_1
\end{array}\right) \hat{x}^{u_3, u_4} \\
& +\left(\begin{array}{c}
	\tilde{P}_1\sum\limits_{i=1}^3 b_i d_i \\
	\tilde{P}_2\sum\limits_{i=1}^3 c_i d_i
\end{array}\right) \hat{u}_3+\left(\begin{array}{c}
\tilde{P}_1\sum\limits_{i=1}^3 b_i e_i \\
\tilde{P}_2\sum\limits_{i=1}^3 c_i e_i
\end{array}\right) \hat{u}_4+\left(\begin{array}{ll}
b_0 & 0 \\
0 & c_0
\end{array}\right)\left(\begin{array}{c}
\hat{\phi}_1 \\
\hat{\phi}_2
\end{array}\right) \\
& +\left(\begin{array}{ll}
b_1 & 0 \\
0 & c_1
\end{array}\right)\left(\begin{array}{l}
\hat{\zeta}_{11} \\
\hat{\zeta}_{21}
\end{array}\right)+\left(\begin{array}{cc}
b_3 & 0 \\
0 & c_3
\end{array}\right)\left(\begin{array}{c}
\hat{\zeta}_{13} \\
\hat{\zeta}_{23}
\end{array}\right)\Bigg].
\end{aligned}
\end{aligned}
$$
Now, we can write the Nash equilibrium point of the followers as
\begin{equation}\label{Nash equilibrium point of the followers}
\begin{aligned}
u_1^*& =-L_{11} \hat{x}^{u_3 , u_4}-L_{12} \hat{u}_3-L_{13} \hat{u}_4+\frac{\tilde{P}_1\left(\sum\limits_{i=1}^3 b_i c_i\right) c_0}{|N|} \hat{\phi}_2 -\frac{\bar{R}_2 b_0}{|N|} \hat{\phi}_1\\
&\quad +\frac{\tilde{P}_1\left( \sum\limits_{i=1}^3 b_i c_i\right) c_1}{|N|} \hat{\zeta}_{21}-\frac{\bar{R}_2 b_1}{|N|} \hat{\zeta}_{11} +\frac{\bar{P}_1\left(\sum\limits_{i=1}^3 b_i c_i\right) c_3}{|N|} \hat{\zeta}_{23}-\frac{\bar{R}_1 b_3}{|N|} \hat{\zeta}_{13}, \\
u_2^*& =-L_{21} \hat{x}^{u_3 , u_4}-L_{22} \hat{u}_3-L_{23} \hat{u}_4+\frac{\tilde{P}_2\left(\sum\limits_{i=1}^3 b_i c_i\right) b_0}{|N|} \hat{\phi}_1 -\frac{\bar{R}_1 c_0}{|N|} \hat{\phi}_2\\
&\quad +\frac{\tilde{P}_2\left( \sum\limits_{i=1}^3 b_i c_i\right) b_1}{|N|} \hat{\zeta}_{11}-\frac{\bar{R}_1 c_1}{|N|} \hat{\zeta}_{21} +\frac{\bar{P}_2\left(\sum\limits_{i=1}^3 b_i c_i\right) b_3}{|N|} \hat{\zeta}_{13}-\frac{\bar{R}_1 c_3}{|N|} \hat{\zeta}_{23},
\end{aligned}
\end{equation}
where
\begin{equation*}
\begin{aligned}
L_{11}&:=\frac{\bar{R}_2 \tilde{P}_1 \bar{B}_1-\tilde{P}_1 \tilde{P}_2 \bar{C}_1\left(\sum\limits_{i=1}^3 b_i c_i\right)}{|N|}, \
L_{12}:=\frac{\bar{R}_2 \tilde{P}_1\left(\sum\limits_{i=1}^3 b_i d_i\right)-\tilde{P}_1 \tilde{P}_2\left(\sum\limits_{i=1}^3 b_i c_i\right)\left(\sum\limits_{i=1}^3 c_i d_i\right)}{|N|}, \\
L_{13}&:=\frac{\bar{R}_2 \tilde{P}_1\left(\sum\limits_{i=1}^3 b_i e_i\right)-\tilde{P}_1 \tilde{P}_2\left(\sum\limits_{i=1}^3 b_i c_i \right)\left(\sum\limits_{i=1}^3 c_i e_i\right)}{|N|}, \
L_{21}:=\frac{\bar{R}_1 \tilde{P}_2 \bar{C}_1-\tilde{P}_1 \tilde{P}_2 \bar{B}_1\left(\sum\limits_{i=1}^3 b_i c_i\right)}{|N|}, \\
L_{22}&:=\frac{\bar{R}_1 \tilde{P}_2\left(\sum\limits_{i=1}^3 c_i d_i\right)-\tilde{P}_1 \tilde{P}_2\left(\sum\limits_{i=1}^3 b_i c_i\right)\left(\sum\limits_{i=1}^3 b_i d_i\right)}{|N|}, \\
L_{23}&:=\frac{\bar{R}_1 \tilde{P}_2\left(\sum\limits_{i=1}^3 c_i e_i\right)-\tilde{P}_1 \tilde{P}_2\left(\sum\limits_{i=1}^3 b_i c_i \right)\left(\sum\limits_{i=1}^3 b_i e_i\right)}{|N|}.
\end{aligned}
\end{equation*}

Next, we need to derive the filtering equations that $(\hat{\phi}_j(\cdot),\hat{\zeta}_{j1}(\cdot),\hat{\zeta}_{j3}(\cdot))$, $j=1,2$ satisfy. Applying Lemma 5.4 of Xiong \cite{Xiong-2008} to (\ref{variphi BSDEs of followers}), we get
\begin{equation}\label{hat variphi BSDEs of followers}
\left\{\begin{array}{l}
d \hat{\phi}_j(t)=\hat{\alpha}_j(t) d t+\hat{\zeta}_{j1}(t) d W_1+\hat{\zeta}_{j3}(t) d W_3, \quad t\in[0,T],\\
\hat{\phi_j}(T)=0, \quad j=1,2,
\end{array}\right.
\end{equation}

Introduce a system of coupled Riccati equations:
\begin{equation}\label{coupled Riccati equations of followers}
\left\{\begin{aligned}
&\dot{\tilde{P}}_1+2 a_0 \tilde{P}_1+\sum\limits_{i=1}^3 \tilde{P}_1 a_i^2+Q_1-\tilde{P}_1 b_0 L_{11}-\tilde{P}_1 c_0 L_{21}\\
&\quad -\left(\sum\limits_{i=1}^3 a_i b_i\right) \tilde{P}_1 L_{11}  -\left(\sum\limits_{i=1}^3 a_i c_i\right) \tilde{P}_1 L_{21}=0,\quad \tilde{P}_1(T)=G_1,  \\
&\dot{\tilde{P}}_2+2 a_0 \tilde{P}_2+\sum\limits_{i=1}^3 \tilde{P}_2 a_i^2+Q_2-\tilde{P}_2 b_0 L_{11}-\tilde{P}_2 c_0 L_{21}\\
&\quad -\left(\sum\limits_{i=1}^3 a_i b_i\right) \tilde{P}_2 L_{11} -\left(\sum\limits_{i=1}^3 a_i c_i\right) \tilde{P}_2 L_{21}=0,\quad \tilde{P}_2(T)=G_2,
\end{aligned}\right.
\end{equation}
where $L_{ji},j=1,2,i=1,2,3$ are defined as above. Since the general solvability of coupled Riccati equations, such as (\ref{coupled Riccati equations of followers}), is difficult, we just discuss its solvability in some special case.

\begin{assumption}\label{A3}
	$R_1(t)=R_2(t)$ and $b_i(t)=c_i(t), i=0,1,2,3$, for $t\in[0,T]$.
\end{assumption}

\newtheorem{lemma}{Lemma}[section]
\begin{lemma}\label{lemma1}
Under Assumption \ref{A3}, (\ref{coupled Riccati equations of followers}) exists a unique solution $(\tilde{P}_1(\cdot),\tilde{P}_2(\cdot))$.
\end{lemma}
The proof of Lemma \ref{lemma1} is left in the Appendix.

According to Lemma \ref{lemma1}, substituting (\ref{filtering relation of followers}), (\ref{filtering comparing of followers}) and (\ref{Nash equilibrium point of the followers}) into (\ref{hat variphi BSDEs of followers}), we obtain
\begin{equation}\label{backward equation 1 of followers}
\left\{\begin{aligned}
d\hat{\phi}_1(t)&= \left(-\tilde{L}_{10} \hat{\phi}_1 -\tilde{L}_{11} \hat{\phi}_2-\tilde{L}_{12} \hat{\zeta}_{11}-\tilde{L}_{13} \hat{\zeta}_{13}-\tilde{L}_{14} \hat{\zeta}_{21}
-\tilde{L}_{15} \hat{\zeta}_{23}-\tilde{L}_{16} \hat{u}_3-\tilde{L}_{17} \hat{u}_4 \right) d t  \\
&\quad +\hat{\zeta}_{11} d W_1+\hat{\zeta}_{13} d W_3,\quad t\in[0,T], \\
\hat{\phi}_1(T) &=0,
\end{aligned}\right.
\end{equation}
where
{\small$$
\begin{aligned}
\tilde{L}_{11}& :=\frac{\left(\sum\limits_{i=1}^3 b_i c_i\right)  \tilde{P}_1^2 b_0 c_0}{|N|}-\frac{ \tilde{P}_1 c_0^2 \bar{R}_1}{|N|}
 +\frac{\left(\sum\limits_{i=1}^3 a_i b_i\right) \left(\sum\limits_{i=1}^3 b_i c_i\right)\tilde{P}_1^2 c_0}{|N|}-\frac{\left(\sum\limits_{i=1}^3 a_i c_i\right) \tilde{P}_1 \bar{R}_1 c_0}{|N|},  \\
\tilde{L}_{10}&:=-\frac{\tilde{P}_1 b_0^2 \bar{R}_2 }{|N|}+\frac{\left(\sum\limits_{i=1}^3 b_i c_i\right) \tilde{P}_1 \tilde{P}_2 c_0 b_0}{|N|}
 -\frac{\left(\sum\limits_{i=1}^3 a_i b_i\right) \tilde{P}_1 \bar{R}_2 b_0}{|N|}+\frac{\left(\sum\limits_{i=1}^3 a_i c_i\right) \left(\sum\limits_{i=1}^3 b_i c_i\right)}{|N|} \tilde{P}_1 \tilde{P}_2 b_0 +a_0, \\
\tilde{L}_{12}& :=-\frac{\tilde{P}_1 b_0 b_1 \bar{R}_2 }{|N|}+\frac{\left(\sum\limits_{i=1}^3 b_i c_i\right) \tilde{P}_1 \tilde{P}_2 c_0 b_1}{|N|}
 -\frac{\left(\sum\limits_{i=1}^3 a_i b_i\right) \tilde{P}_1 \bar{R}_2 b_1}{|N|}+\frac{\left(\sum\limits_{i=1}^3 a_i c_i\right) \left(\sum\limits_{i=1}^3 b_i c_i\right)}{|N|} \tilde{P}_1 \tilde{P}_2 b_1 +a_1, \\
\end{aligned}
$$
$$
\begin{aligned}
\tilde{L}_{13}& :=-\frac{\tilde{P}_1 b_0 b_3 \bar{R}_2 }{|N|}+\frac{\left(\sum\limits_{i=1}^3 b_i c_i\right) \tilde{P}_1 \tilde{P}_2 c_0 b_3}{|N|}
 -\frac{\left(\sum\limits_{i=1}^3 a_i b_i\right) \tilde{P}_1 \bar{R}_2 b_3}{|N|}+\frac{\left(\sum\limits_{i=1}^3 a_i c_i\right) \left(\sum\limits_{i=1}^3 b_i c_i\right)}{|N|} \tilde{P}_1 \tilde{P}_2 b_3 +a_3, \\
\tilde{L}_{14}& :=\frac{\left(\sum\limits_{i=1}^3 b_i c_i\right)  \tilde{P}_1^2 b_0 c_1}{|N|}-\frac{ \tilde{P}_1 c_0 c_1 \bar{R}_1}{|N|}
 +\frac{\left(\sum\limits_{i=1}^3 a_i b_i\right) \left(\sum\limits_{i=1}^3 b_i c_i\right)\tilde{P}_1^2 c_1}{|N|}-\frac{\left(\sum\limits_{i=1}^3 a_i c_i\right) \tilde{P}_1 \bar{R}_1 c_1}{|N|},   \\
\tilde{L}_{15}& :=\frac{\left(\sum\limits_{i=1}^3 b_i c_i\right)  \tilde{P}_1^2 b_0 c_3}{|N|}-\frac{ \tilde{P}_1 c_0 c_3 \bar{R}_1}{|N|}
 +\frac{\left(\sum\limits_{i=1}^3 a_i b_i\right) \left(\sum\limits_{i=1}^3 b_i c_i\right)\tilde{P}_1^2 c_3}{|N|}-\frac{\left(\sum\limits_{i=1}^3 a_i c_i\right) \tilde{P}_1 \bar{R}_1 c_3}{|N|},  \\
\tilde{L}_{16}& :=-\tilde{P}_1 b_0 L_{12}-\tilde{P}_1 c_0 L_{22}+d_0 \tilde{P}_1-\left(\sum_{i 1}^3 a_i b_i\right) \tilde{P}_1 L_{12}
 -\left(\sum_{i 1}^3 a_i c_i\right) \tilde{P}_1 L_{22}+\left(\sum_{i=1}^3 a_i d_i\right) \tilde{P}_1, \\
\tilde{L}_{17}& :=-\tilde{P}_1 b_0 L_{13}-\tilde{P}_1 c_0 L_{23}+e_0 \tilde{P}_1-\left(\sum_{i=1}^3 a_i b_i\right) \tilde{P}_1 L_{13}
 -\left(\sum_{i=1}^3 a_i c_i\right) \tilde{P}_1 L_{23}+\left(\sum_{i=1}^3 a_i e_i\right) \tilde{P}_1.
\end{aligned}
$$}
and
\begin{equation}\label{backward equation 2 of followers}
\left\{\begin{aligned}
d\hat{\phi}_2(t)&= \left(-\tilde{L}_{20} \hat{\phi}_2 -\tilde{L}_{21} \hat{\phi}_1-\tilde{L}_{22} \hat{\zeta}_{11}-\tilde{L}_{23} \hat{\zeta}_{13}-\tilde{L}_{24} \hat{\zeta}_{21}-\tilde{L}_{25} \hat{\zeta}_{23}-\tilde{L}_{26} \hat{u}_3-\tilde{L}_{27} \hat{u}_4 \right) d t  \\
&\quad +\hat{\zeta}_{21} d W_1+\hat{\zeta}_{23} d W_3,\quad t\in[0,T], \\
\hat{\phi}_2(T) &= 0,
\end{aligned}\right.
\end{equation}
where
{\small$$
\begin{aligned}
\tilde{L}_{20}& :=\frac{\left(\sum\limits_{i=1}^3 b_i c_i\right)  \tilde{P}_1 \tilde{P}_2 b_0 c_0}{|N|}-\frac{ \tilde{P}_2 c_0^2 \bar{R}_1}{|N|}
 +\frac{\left(\sum\limits_{i=1}^3 a_i b_i\right) \left(\sum\limits_{i=1}^3 b_i c_i\right)\tilde{P}_1 \tilde{P}_2 c_0}{|N|}-\frac{\left(\sum\limits_{i=1}^3 a_i c_i\right) \tilde{P}_2 \bar{R}_1 c_0}{|N|} +a_0, \\
\tilde{L}_{21}& :=-\frac{\tilde{P}_2 b_0 \bar{R}_2 b_0}{|N|}+\frac{\left(\sum\limits_{i=1}^3 b_i c_i\right) \tilde{P}_2^2 c_0 b_0}{|N|}
 -\frac{\left(\sum\limits_{i=1}^3 a_i b_i\right) \tilde{P}_2 \bar{R}_2 b_0}{|N|}+\frac{\left(\sum\limits_{i=1}^3 a_i c_i\right) \left(\sum\limits_{i=1}^3 b_i c_i\right)}{|N|} \tilde{P}_2^2 b_0, \\
\tilde{L}_{22}& :=-\frac{\tilde{P}_2 b_0 \bar{R}_2 b_1}{|N|}+\frac{\left(\sum\limits_{i=1}^3 b_i c_i\right) \tilde{P}_2^2 c_0 b_1}{|N|}
 -\frac{\left(\sum\limits_{i=1}^3 a_i b_i\right) \tilde{P}_2 \bar{R}_2 b_1}{|N|}+\frac{\left(\sum\limits_{i=1}^3 a_i c_i\right)\left(\sum\limits_{i=1}^3 b_i c_i\right)}{|N|} \tilde{P}_2^2 b_1, \\
\tilde{L}_{23}& :=-\frac{\tilde{P}_2 b_0 \bar{R}_2 b_3}{|N|}+\frac{\left(\sum\limits_{i=1}^3 b_i c_i\right) \tilde{P}_2^2 c_0 b_3}{|N|}
 -\frac{\left(\sum\limits_{i=1}^3 a_i b_i\right) \tilde{P}_2 \bar{R}_2 b_3}{|N|}+\frac{\left(\sum\limits_{i=1}^3 a_i c_i\right)\left(\sum\limits_{i=1}^3 b_i c_i\right)}{|N|} \tilde{P}_2^2 b_3, \\
\end{aligned}
$$
$$
\begin{aligned}
\tilde{L}_{24}& :=\frac{\left(\sum\limits_{i=1}^3 b_i c_i\right)  \tilde{P}_1 \tilde{P}_2 b_0 c_1}{|N|}-\frac{ \tilde{P}_2 c_0 c_1 \bar{R}_1}{|N|}
 +\frac{\left(\sum\limits_{i=1}^3 a_i b_i\right) \left(\sum\limits_{i=1}^3 b_i c_i\right)\tilde{P}_1 \tilde{P}_2 c_1}{|N|}-\frac{\left(\sum\limits_{i=1}^3 a_i c_i\right) \tilde{P}_2 \bar{R}_1 c_1}{|N|} +a_1, \\
\tilde{L}_{25}& :=\frac{\left(\sum\limits_{i=1}^3 b_i c_i\right)  \tilde{P}_1 \tilde{P}_2 b_0 c_3}{|N|}-\frac{ \tilde{P}_2 c_0 c_3 \bar{R}_1}{|N|}
 +\frac{\left(\sum\limits_{i=1}^3 a_i b_i\right) \left(\sum\limits_{i=1}^3 b_i c_i\right)\tilde{P}_1 \tilde{P}_2 c_3}{|N|}-\frac{\left(\sum\limits_{i=1}^3 a_i c_i\right) \tilde{P}_2 \bar{R}_1 c_3}{|N|} +a_3, \\
\tilde{L}_{26}& :=-\tilde{P}_2 b_0 L_{12}-\tilde{P}_2 c_0 L_{22}+d_0 \tilde{P}_2-\left(\sum_{i 1}^3 a_i b_i\right) \tilde{P}_2 L_{12}
 -\left(\sum_{i 1}^3 a_i c_i\right) \tilde{P}_2 L_{22}+\left(\sum_{i=1}^3 a_i d_i\right) \tilde{P}_2, \\
\tilde{L}_{27}& :=-\tilde{P}_2 b_0 L_{13}-\tilde{P}_2 c_0 L_{23}+e_0 \tilde{P}_2-\left(\sum_{i=1}^3 a_i b_i\right) \tilde{P}_2 L_{13}
 -\left(\sum_{i=1}^3 a_i c_i\right) \tilde{P}_2 L_{23}+\left(\sum_{i=1}^3 a_i e_i\right) \tilde{P}_2.
\end{aligned}
$$}

Noticing that (\ref{backward equation 1 of followers}) and (\ref{backward equation 2 of followers}) are coupled BSDEs. For the solvability of them, we set
\begin{equation}\label{Psi definition of folowers}
\begin{aligned}
\Phi:=\left(\begin{array}{l}
\phi_1 \\
\phi_2
\end{array}\right),\quad \zeta_1:=\left(\begin{array}{l}
\zeta_{11} \\
\zeta_{13}
\end{array}\right),\quad \zeta_3:=\left(\begin{array}{l}
\zeta_{21} \\
\zeta_{23}
\end{array}\right),
\end{aligned}
\end{equation}
and thus
\begin{equation}\label{high dimension backward equation of followers}
\left\{\begin{aligned}
d \hat{\Phi}(t)&= -\left[\left(\begin{array}{ll}
\tilde{L}_{10} & \tilde{L}_{11} \\
\tilde{L}_{21} & \tilde{L}_{20}
\end{array}\right) \hat{\Phi}+\left(\begin{array}{ll}
\tilde{L}_{12} & \tilde{L}_{13} \\
\tilde{L}_{22} & \tilde{L}_{23}
\end{array}\right) \hat{\zeta}_1
+\left(\begin{array}{ll}
\tilde{L}_{14} & \tilde{L}_{15} \\
\tilde{L}_{24} & \tilde{L}_{26}
\end{array}\right) \hat{\zeta}_2  \right. \\
&\qquad \left. +\left(\begin{array}{l}
\tilde{L}_{16} \\
\tilde{L}_{26}
\end{array}\right) \hat{u}_3+\left(\begin{array}{l}
\tilde{L}_{17} \\
\tilde{L}_{27}
\end{array}\right) \hat{u}_4\right] d t
+\hat{\zeta}_1 d W_1+\hat{\zeta}_3 d W_3,\quad t\in[0,T], \\
\hat{\Phi}(T)&= O_{2 \times 1}.
\end{aligned}\right.
\end{equation}
This is a standard BSDE, which admits a unique $\mathcal{G}_t^1$-adapted solution triple $(\hat{\Phi}(\cdot), \hat{\zeta_1}(\cdot), \hat{\zeta}_3(\cdot))$.

Putting (\ref{Nash equilibrium point of the followers}) into the forward equation of (\ref{optimality system of followers}), and applying Lemma 5.4 of \cite{Xiong-2008}, we get
\begin{equation}\label{filtering state of followers}
\left\{\begin{aligned}
d\hat{x}^{u_3,u_4}(t) &= \Delta(t;\hat{u}_3,\hat{u}_4) d t +\sum\limits_{i=1,3}\Upsilon_i(t;\hat{u}_3,\hat{u}_4)dW_i, \quad t\in[0,T],\\
\hat{x}^{u_3,u_4}(0) &= x_0,
\end{aligned}\right.
\end{equation}
where
$$
\begin{aligned}
&\Delta(t;\hat{u}_3,\hat{u}_4):=\left(a_0-b_0 L_{11}-c_i L_{21}\right) \hat{x}^{u_3, u_4} +\frac{-\bar{R}_2 b_0^2+\left(\sum\limits_{i=1}^3 b_i c_i \right) \tilde{P}_2 b_0 c_0}{|N|} \hat{\phi}_1 \\
&\qquad +\frac{-\bar{R}_1 c_0^2+\left(\sum\limits_{i=1}^3 b_i c_i\right) \tilde{P}_1 b_0 c_0}{|N|} \hat{\phi}_2+\frac{-\bar{R}_2 b_1 b_0+\left(\sum\limits_{i=1}^3 b_i c_i\right) \tilde{P}_2 b_1 c_0}{|N|} \hat{\zeta}_{11} \\
&\qquad +\frac{-\bar{R}_2 b_3 b_0+\left(\sum\limits_{i=1}^3 b_i c_i\right) \tilde{P}_2 b_3 c_0}{|N|} \hat{\zeta}_{13}+\frac{-\bar{R}_1 c_0 c_1+\left(\sum\limits_{i=1}^3 b_i c_i\right) \tilde{P}_1 c_1 b_0}{|N|} \hat{\zeta}_{21} \\
\end{aligned}
$$
$$
\begin{aligned}
&\qquad +\frac{-\bar{R}_1 c_3 c_0+\left(\sum\limits_{i=1}^3 b_i c_i\right) \tilde{P}_1 c_3 b_0}{|N|} \hat{\zeta}_{23}+\left(d_0-b_0 L_{12}-c_0 L_{22}\right) \hat{u}_3 +\left(e_0-b_0 L_{13}-c_0 L_{23}\right) \hat{u}_4,\\
&\Upsilon_i(t;\hat{u}_3,\hat{u}_4):=\left(a_i-b_i L_{11}-c_i L_{21}\right) \hat{x}^{u_3, u_4} +\frac{-\bar{R}_2 b_0 b_i+\left(\sum\limits_{i=1}^3 b_i c_i \right) \tilde{P}_2 b_0 c_i}{|N|} \hat{\phi}_1 \\
&\qquad +\frac{-\bar{R}_1 c_0 c_i+\left(\sum\limits_{i=1}^3 b_i c_i\right) \tilde{P}_1 c_0 b_i}{|N|} \hat{\phi}_2+\frac{-\bar{R}_2 b_1 b_i+\left(\sum\limits_{i=1}^3 b_i c_i\right) \tilde{P}_2 b_1 c_i}{|N|} \hat{\zeta}_{11} \\
&\qquad +\frac{-\bar{R}_2 b_3 b_i+\left(\sum\limits_{i=1}^3 b_i c_i\right) \tilde{P}_2 b_3 c_i}{|N|} \hat{\zeta}_{13}+\frac{-\bar{R}_1 c_1 c_i+\left(\sum\limits_{i=1}^3 b_i c_i\right) \tilde{P}_1 c_1 b_i}{|N|} \hat{\zeta}_{21} \\
&\qquad +\frac{-\bar{R}_1 c_3 c_i+\left(\sum\limits_{i=1}^3 b_i c_i\right) \tilde{P}_1 c_3 b_i}{|N|} \hat{\zeta}_{23}+\left(d_i-b_i L_{12}-c_i L_{22}\right) \hat{u}_3+\left(e_i-b_i L_{13}-c_i L_{23}\right) \hat{u}_4.
\end{aligned}
$$
It is obvious that (\ref{filtering state of followers}) admits a unique $\mathcal{G}_t^1$-adapted solution $\hat{x}^{u_3, u_4}(\cdot)$.

We summarize the above statements in the following theorem.
\begin{thm}\label{thm2}
Let Assumption \ref{A1}, Assumption \ref{A2} hold and $P_j(\cdot),j=1,2$ satisfy (\ref{coupled Riccati equations of followers}). For chosen $( u_3(\cdot), u_4(\cdot) )\in \mathcal{U}_3 \times \mathcal{U}_4$ of the leaders, let $\left(u_1^*(\cdot), u_2^*(\cdot)\right)$ be a Nash equilibrium point of the followers, then it has the state-estimate feedback representation of (\ref{Nash equilibrium point of the followers}), where $(\hat{x}^{u_3, u_4}(\cdot), \hat{\phi}_j(\cdot), \hat{\zeta}_{j1}(\cdot), \hat{\zeta}_{j3}(\cdot))$ is determined by (\ref{high dimension backward equation of followers}) and (\ref{filtering state of followers}), for $j=1,2$.
\end{thm}

\begin{Remark}\label{remark1}
For any given $( u_3(\cdot), u_4(\cdot) )\in \mathcal{U}_3 \times \mathcal{U}_4$, since $(\hat{\Phi}(\cdot), \hat{\zeta}_1(\cdot), \hat{\zeta}_3(\cdot))$, $\hat{x}^{u_3, u_4}(\cdot)$ are unique solutions to (\ref{high dimension backward equation of followers}), (\ref{filtering state of followers}) respectively, and $\tilde{P}_j(\cdot), j=1,2$ are unique solutions to (\ref{coupled Riccati equations of followers}) under Assumption \ref{A3}, the Nash equilibrium point $\left(u_1^*(\cdot), u_2^*(\cdot)\right)$ of the followers' problem is unique by combining with the sufficient condition in Theorem \ref{thm1}.
\end{Remark}

\subsection{Problem of the leaders}

After knowing that the followers would take their Nash equilibrium point (\ref{Nash equilibrium point of the followers}), the leaders seek a pair $\left(u_3^*(\cdot), u_4^*(\cdot)\right)$ satisfying (\ref{Nash game of leaders}). Now, the leaders' ``state" equation can be written as
\begin{equation}\label{state of leaders}
\left\{\begin{aligned}
d x^{u_3, u_4}(t)&= \Big(a_0 x^{u_3, u_4}+M_{01} \hat{x}^{u_3, u_4}+N_{00}^T \hat{\Phi} +N_{01}^T \hat{\zeta}_1 +N_{03}^T \hat{\zeta}_3+d_0 u_3+M_{04} \hat{u}_3\\
&\quad +e_0 u_4+M_{05} \hat{u}_4\Big) d t +\sum\limits_{i=1}^3\Big( a_i x^{u_3 u_4}+M_{i 1} \hat{x}^{u_3 u_4}+N_{i0}^T \hat{\Phi} +N_{i1}^T \hat{\zeta}_1 \\
&\quad +N_{i3}^T \hat{\zeta}_3+d_i u_3+M_{i 4} \hat{u}_3+e_i u_4+M_{i 5} \hat{u}_4\Big) d W_i,  \\
d \hat{\Phi}(t)&= -\big( A_0 \hat{\Phi}+ A_1 \hat{\zeta}_1+A_3 \hat{\zeta}_3+f_1 \hat{u}_3  +f_2 \hat{u}_4\big) d t +\hat{\zeta}_1 d W_1+\hat{\zeta}_3 d W_3, \quad t\in[0,T],\\
x^{u_3, u_4}(0)& =x_0, \quad \hat{\Phi}(T)= O_{2 \times 1},
\end{aligned}\right.  \\
\end{equation}
where for $i=0,1,2,3$,
$$
\begin{aligned}
M_{i 1}& :=-b_i L_{11}-c_i L_{21},\quad  M_{i 4}  =-b_i L_{12}-c_i L_{22}, \quad  M_{i 5} = -b_i L_{13}-c_i L_{23}, \\
M_{i 2}& :=\frac{-\bar{R}_2 b_i+\left(\sum\limits_{l=1}^3 b_l c_l\right) \tilde{P}_2 c_i}{|N|},\quad M_{i 3} :=\frac{-\bar{R}_1 c_i+\left(\sum\limits_{l=1}^3 b_l c_l\right) \tilde{P}_1 b_i}{|N|},\\
N_{i l}& :=\left(\begin{array}{c}
           M_{i 2}b_l \\
           M_{i 3}c_l
          \end{array}\right),\ l=0,1,3,\quad
A_0 :=\left(\begin{array}{ll}
\tilde{L}_{10} & \tilde{L}_{11} \\
\tilde{L}_{21} & \tilde{L}_{20}
\end{array}\right),\quad
A_1 :=\left(\begin{array}{ll}
\tilde{L}_{12} & \tilde{L}_{13} \\
\tilde{L}_{22} & \tilde{L}_{23}
\end{array}\right),\\
A_3& :=\left(\begin{array}{ll}
\tilde{L}_{14} & \tilde{L}_{15} \\
\tilde{L}_{24} & \tilde{L}_{25}
\end{array}\right), \quad
f_1 :=\left(\begin{array}{c}
\tilde{L}_{16} \\
\tilde{L}_{26}
\end{array}\right),\quad f_2:=\left(\begin{array}{l}
\tilde{L}_{17} \\
\tilde{L}_{27}
\end{array}\right).
\end{aligned}
$$
Note that (\ref{state of leaders}) is a CMF-FBSDE, whose uniquely solvability of its $\mathcal{F}_t$-adapted solution quadruple $(x^{u_3,u_4}(\cdot),\hat{\Phi}(\cdot), \hat{\zeta}_1(\cdot), \hat{\zeta}_3(\cdot))$ can be similarly guaranteed as in the previous subsection.

The cost functionals of the leaders could be written as, for $j=3,4$,
\begin{equation}\label{cost of leaders}
\hat{J}_j(u_3(\cdot), u_4(\cdot)):= \frac{1}{2} \mathbb{E}\left[\int_0^T\left[ Q_j(t) x^{u_3, u_4}(t)^2+ R_j(t) u_j(t)^2\right] d t
+ G_j x^{u_3, u_4}(T)^2\right].
\end{equation}

For the leaders' problem, we introduce the following assumption.
\begin{assumption}\label{A4}
$Q_j(t) \geq 0, R_j(t) > 0, \forall t \in[0, T]$, and $G_j \geq 0$, for $j=3,4$.
\end{assumption}

The two leaders encounter a partial information LQ nonzero-sum stochastic differential Nash game. We first have the following result whose proof is also left to the Appendix.
\begin{thm}\label{thm3}
Let Assumptions \ref{A1}, \ref{A2}, \ref{A4} hold. A pair $\left(u_3^*(\cdot), u_4^*(\cdot)\right)$ is a Nash equilibrium point of the leaders, if and only if
\begin{equation}\label{SMP of leaders}
\begin{aligned}
& u_3^*=-R_3^{-1}\bigg[d_0 \check{z}_3+M_{04} \check{\hat{z}}_3+\sum_{i=1}^3\left(d_i \check{q}_{3 i}+M_{i 4} \check{\hat{q}}_{3 i}\right)+f_1^{\top} \check{\hat{y}}_3\bigg], \\
& u_4^*=-R_4^{-1}\bigg[e_0 \check{z}_4+M_{05} \check{\hat{z}}_4+\sum_{j=1}^3\left(e_i \check{q}_{4 i}+M_{i 5} \check{\hat{q}}_{4 i}\right)+f_2^{\top} \check{\hat{y}}_4\bigg], \\
& \hat{u}_3^*=-R_3^{-1}\bigg[\left(d_0+M_{0 4}\right)\check{\hat{z}}_3+\sum_{i=1}^3\left(d_i+M_{i 4}\right) \check{\hat{q}}_{3 i}+f_1^{\top} \check{\hat{y}}_3\bigg], \\
& \hat{u}_4^*=-R_4^{-1}\bigg[\left(e_0+M_{05}\right) \check{\hat{z}}_4+\sum_{i=1}^3\left(e_i+M_{i 5}\right) \hat{q}_{4 i}+f_2^{\top} \check{\hat{y}}_4\bigg],
\end{aligned}
\end{equation}
where the $\mathcal{F}_t$-adapted process quadruple $\left(y_j(\cdot), z_j(\cdot), q_{j 1}(\cdot), q_{j 2}(\cdot), q_{j 3}(\cdot)\right) \in \mathbb{R}^2 \times \mathbb{R} \times \mathbb{R} \times \mathbb{R} \times \mathbb{R}$ satisfies the adjoint CMF-FBSDEs:
\begin{equation}\label{adjoint CMF-FBSDE of leaders}
\left\{\begin{aligned}
d y_j(t) &= \bigg(A_0^\top y_j+N_{00} z_j+\sum_{i=1}^3 N_{i 0} q_{j i}\bigg) d t +\bigg(A_1^\top y_j+N_{01} z_j+\sum_{i=1}^3 N_{i 1} q_{j i}\bigg) d W_1 \\
&\quad +\bigg(A_3^\top y_j+N_{03} z_j+\sum_{i=1}^3 N_{i 3} q_{j i}\bigg) d W_3,  \\
d z_j(t) &= \big(-a_0 z_j-a_1 q_{j1}-a_2 q_{j2}-a_3 q_{j3}-M_{01} \hat{z}_j-M_{11} \hat{q}_{j1} \\
&\quad -M_{21} \hat{q}_{j2} - M_{31} \hat{q}_{j3}-Q_j x^{*}\big)dt+\sum\limits_{i=1}^3 q_{ji} d W_i,\quad t\in[0,T], \\
 y_j(0)&=0, \quad z_j(T)=G_j x^*(T).
\end{aligned}\right.
\end{equation}
In the above, we have denoted $x^*(\cdot)\equiv x^{u_3^*,u_4^*}(\cdot)$.
\end{thm}

Next, we want to obtain the feedback form of $\left(u_3^*(\cdot), u_4^*(\cdot)\right)$. However, the situation which the leaders encounter is different from that of the followers, since now the ``states" of the leaders are the process quadruple $(x^{u_3, u_4}(\cdot), \hat{\Phi}(\cdot), \hat{\zeta}_1(\cdot), \hat{\zeta}_3(\cdot))$ which satisfy the CMF-FBSDE (\ref{state of leaders}), not (\ref{state}) for just the initial state process $x^{u_3, u_4}(\cdot)$! To overcome this difficulty, inspired by \cite{Yong-2002}, let
\begin{equation}\label{X,Y,Z}
\begin{aligned}
& X:=\left(\begin{array}{l}
x^{u_3,u_4} \\
y_3 \\
y_4
\end{array}\right),\ Y:=\left(\begin{array}{l}
\hat{\Phi} \\
z_3 \\
z_4
\end{array}\right),\ Z_1:=\left(\begin{array}{l}
\hat{\zeta}_1 \\
q_{31} \\
q_{41}
\end{array}\right),\ Z_2:=\left(\begin{array}{l}
0_{2 \times 1} \\
q_{32} \\
q_{42}
\end{array}\right),\ Z_3=\left(\begin{array}{l}
\hat{\zeta}_3 \\
q_{33} \\
q_{43}
\end{array}\right),
\end{aligned}
\end{equation}
where $\left(X(\cdot), Y(\cdot), Z_1(\cdot), Z_2(\cdot) , Z_3(\cdot) \right)\in \mathbb{R}^5 \times \mathbb{R}^4 \times \mathbb{R}^4 \times \mathbb{R}^4 \times \mathbb{R}^4$,
then (\ref{state of leaders}) and (\ref{adjoint CMF-FBSDE of leaders}) can be rewritten as
\begin{equation}\label{Optimality system of leaders}
\left\{ \begin{aligned}
dX(t) & =\left(\mathcal{A}_0 X+\mathcal{M}_{01} \hat{X}+\mathcal{H}_0 \check{\hat{X}}+\mathcal{N}_{00} Y+\mathcal{B}_{00} \hat{Y}+\mathcal{C}_{00} \check{\hat{Y}} +\mathcal{N}_{01}Z_1 +\mathcal{B}_{01} \check{Z}_1+\mathcal{C}_{01} \check{\hat{Z}}_1 \right.\\
&\quad \left.+\mathcal{N}_{02} Z_2+\mathcal{B}_{02} \check{Z}_2+\mathcal{C}_{02} \check{\hat{Z}}_2+\mathcal{N}_{03} Z_3 +\mathcal{B}_{03} \check{Z}_3+ \mathcal{C}_0^3 \check{\hat{Z}}_3 \right)dt \\
&\quad +\sum_{i=1}^3 \left(\mathcal{A}_i X + \mathcal{B}_{i 2} \check{Z}_2 +\mathcal{H}_i \check{\hat{X}}+\mathcal{N}_{i 0} Y+\mathcal{B}_{i 0} \check{Y}+\mathcal{C}_{i0} \check{\hat{Y}}+\mathcal{N}_{i 1} Z_1+\mathcal{B}_{i1}\check{Z}_1\right. \\
&\quad \left. + \mathcal{C}_{i1} \check{\hat{Z}}_1 +\mathcal{N}_{i 2} Z_2 +\mathcal{B}_{i2}\check{Z}_2+\mathcal{C}_{i2} \check{\hat{Z}}_2+\mathcal{N}_{i 3} Z_3 +\mathcal{B}_{i3}\check{Z}_3+ \mathcal{C}_{i3} \check{\hat{Z}}_3\right)dW_i, \\
d Y(t) & =-\left(\mathcal{Q} X+\check{\hat{\mathcal{Q}}} \check{\hat{X}}+\overline{\mathcal{A}}_0 Y+\overline{\mathcal{M}}_{01} \hat{Y}+\mathcal{I}_0 \check{\hat{Y}} +\overline{\mathcal{A}}_1 Z_1 +\overline{\mathcal{M}}_{11} \hat{Z}_1 +\mathcal{I}_1 \check{\hat{Z}}_1  \right. \\
&\qquad \left. +\overline{\mathcal{A}}_2 Z_2 +\overline{\mathcal{M}}_{21} \hat{Z}_2+\mathcal{ I}_2 \check{\hat{Z}}_2+\overline{A}_3 Z_3+\overline{M}_{31} \hat{Z}_3+\mathcal{I}_3 \check{\hat{Z}}_3 \right)dt
+\sum_{i=1}^3 Z_i d W_i, \quad t\in[0,T],\\
X(0)& =\mathcal{X}_0, \quad Y(T)=\mathcal{G} X(T),
\end{aligned}\right.
\end{equation}
where, for $i=0,1,3$, we have denoted
{\small$$
\begin{aligned}
&\mathcal{A}_i:=\left(\begin{array}{lll}
a_i & & \\
& A_i^\top & \\
& & A_i^{\top}
\end{array}\right),\ \mathcal{N}_{i 0}:=\left(\begin{array}{lll}
N_{i0}^\top & &\\
& N_{0i} &\\
& & N_{0i}
\end{array}\right),\ \mathcal{N}_{i 1}:=\left(\begin{array}{lll}
N_{i1}^\top & &\\
& N_{1i} &\\
& & N_{1i}
\end{array}\right),\\
&\mathcal{N}_{i 2}:=\left(\begin{array}{lll}
O_{1 \times 2} & &\\
& N_{2i} &\\
& & N_{2i}
\end{array}\right),\ \mathcal{N}_{i 3}:=\left(\begin{array}{lll}
N_{i3}^\top & &\\
& N_{3i} &\\
& & N_{3i}
\end{array}\right),\ \mathcal{N}_{2 i}:=\left(\begin{array}{ll}
N_{2i}^\top & \\
& O_{4 \times 2}
\end{array}\right), \\
&\overline{\mathcal{A}}_i:=\left(\begin{array}{lll}
A_i & & \\
& a_i & \\
& & a_i
\end{array}\right);
\end{aligned}
$$}
for $i=0,1,2,3,4$, we have denoted
{\small$$
\begin{aligned}
&\mathcal{M}_{i 1}:=\left(\begin{array}{ll}
M_{i 1} & \\
& O_{4 \times 4}
\end{array}\right),\ \mathcal{H}_i:=\left(\begin{array}{lll}
& {-\left(d_i+M_{i4}\right)R_3^{-1}f_1^\top} & {-\left(e_i+M_{i5}\right)R_4^{-1}+f_2^\top} \\
O_{4 \times 1} & & \\
\end{array}\right),\\
&\overline{\mathcal{M}}_{i 1}=\left(\begin{array}{ccc}
O_{2 \times 2} & & \\
& M_{i 1} & \\
& & M_{i 1}
\end{array}\right),\ \mathcal{I}_i:=\left(\begin{array}{lll}
& -\left(d_i+M_{i 4}\right) R_3^{-1} f_1 & -\left(e_i+M_{i j}\right) R_4^{-1} f_2 \\
O_{2 \times 2} & &
\end{array}\right),\\
\end{aligned}
$$
$$
\begin{aligned}
& \overline{\mathcal{D}}_i:=\left(\begin{array}{c}
O_{2 \times 1} \\
d_i \\
0
\end{array}\right),\ \overline{\mathcal{M}}_{i 4}:=\left(\begin{array}{c}
O_{2 \times 1} \\
M_{i 4} \\
0
\end{array}\right),\ \overline{\mathcal{E}}_i:=\left(\begin{array}{c}
O_{2 \times 1} \\
0 \\
e_i
\end{array}\right),\ \overline{\mathcal{M}}_{i 5}:=\left(\begin{array}{c}
O_{2 \times 1} \\
0 \\
M_{i 5}
\end{array}\right);
\end{aligned}
$$}
for $i=0,1,2,3$, $j=0,1,2,3$, we have denoted
{\small$$
\begin{aligned}
& \mathcal{B}_{ij}:=\left(\begin{array}{lll}
& {-d_i d_j R_3^{-1}} & {-e_i e_j R_4^{-1}} \\
O_{4 \times 2} & & \\
\end{array}\right),\\
& \mathcal{C}_{i j}:=\left(\begin{array}{lll}
& -\left(d_i M_{j 4}+M_{i 4} d_j+M_{i 4} M_{j 4}\right) R_3^{-1} & -\left(e_i M_{j 5}+M_{i 5} e_j+M_{i 5} M_{j 5}\right) R_4^{-1} \\
O_{4 \times 2} & &
\end{array}\right);
\end{aligned}
$$}
and in addition,
{\small$$
\begin{aligned}
& \mathcal{N}_{22}:=O_{5 \times 4},\ \mathcal{Q}:=\left(\begin{array}{ll}
& O_{2 \times 4} \\
Q_3 &  \\
Q_4 &
\end{array}\right),\ \mathcal{F}_1:=\left(\begin{array}{c}
0 \\
f_1 \\
O_{2 \times 1}
\end{array}\right),\\
& \mathcal{F}_2:=\left(\begin{array}{c}
O_{3 \times 1} \\
f_2
\end{array}\right),\ \mathcal{X}_0:=\left(\begin{array}{l}
x_0 \\
O_{4 \times 1}
\end{array}\right),\ G:=\left(\begin{array}{ll}
& O_{2 \times 4} \\
G_3 &  \\
G_4 &
\end{array}\right).
\end{aligned}
$$}
Using above notations, (\ref{SMP of leaders}) also can be written as
\begin{equation}\label{SMP of leaders-}
\begin{aligned}
& u_3^*=-R_3^{-1}\bigg[\overline{\mathcal{D}}_0^\top \check{Y}_t+\overline{\mathcal{M}}_{04}^\top \check{\hat{Y}}_t+\sum_{j=1}^3\left(\overline{\mathcal{D}}_j^\top \check{Z}_{t}^j+\overline{\mathcal{M}}_{j 4}^\top \check{\hat{Z}}_{t}^j \right)+\mathcal{F}_1^\top \check{\hat{X}}_t\bigg], \\
& u_4^*=-R_4^{-1}\bigg[\overline{\mathcal{E}}_0^\top \check{Y}_t+\overline{\mathcal{M}}_{05}^\top \check{\hat{Y}}_t+\sum_{j=1}^3\left(\overline{\mathcal{E}}_j^\top \check{Z}_{t}^{j}+\overline{\mathcal{M}}_{j5}^\top \check{\hat{Z}}_{t}^j \right)+\mathcal{F}_2^\top \check{\hat{X}}_t\bigg].
\end{aligned}
\end{equation}
Applying $\mathbb{E}\left[\cdot \mid \mathcal{G}_t^1\right]$, $\mathbb{E}\left[\cdot \mid \mathcal{G}_t^2\right]$ and $\mathbb{E}\left[\mathbb{E} \left[ \cdot \mid \mathcal{G}_t^1\right] \mid \mathcal{G}_t^2 \right]$ to (\ref{Optimality system of leaders}), we can obtain
\begin{equation}\label{hat X}
\left\{\begin{aligned}
d \hat{X}(t)&=\left[\left(\mathcal{A}_0+\mathcal{M}_{01}\right) \hat{X}+\mathcal{H}_0 \check{\hat{X}}+\mathcal{N}_{00} \hat{Y}+\left(\mathcal{B}_{00}+\mathcal{C}_{00}\right) \check{\hat{Y}}+\mathcal{N}_{01} \hat{Z}_1+ \left(\mathcal{B}_{01}+\mathcal{C}_{01}\right) \check{\hat{Z}}_1\right. \\
&\quad \left.+\mathcal{N}_{02} \hat{Z}_2+\left(\mathcal{B}_{02}+\mathcal{C}_{02}\right) \check{\hat{Z}}_2+\mathcal{N}_{03} \hat{Z}_3+\left(\mathcal{B}_{03}+\mathcal{C}_{03}\right) \check{\hat{Z}}_3\right]dt \\
&\quad +\sum_{i=1,3}\left[\left(\mathcal{A}_i+\mathcal{M}_{i 1}\right) \hat{X}+\mathcal{H}_i \check{\hat{X}}+\mathcal{N}_{i0} \hat{Y}+\left(\mathcal{B}_{i 0}+\mathcal{C}_{i 0}\right) \check{\hat{Y}}+\mathcal{N}_{i1} \hat{Z}_1\right. \\
&\quad \left. +\left(\mathcal{B}_{i1}+\mathcal{C}_{i1}\right) \check{\hat{Z}}_1+\mathcal{N}_{i 2} \hat{Z}_2+\left(\mathcal{B}_{i 2}+\mathcal{C}_{i 2}\right) \check{\hat{Z}}_2+\mathcal{N}_{i 3} \hat{Z}_3
+\left(\mathcal{B}_{i 3}+\mathcal{C}_{i 3}\right) \check{\hat{Z}}_3 \right]dW_i, \ t\in[0,T],\\
\hat{X}(0) &=\mathcal{X}_0,
\end{aligned}\right.
\end{equation}
\begin{equation}\label{check X}
\left\{\begin{aligned}
d \check{X}(t) &=\left[\mathcal{A}_0 \check{X}+\left(\mathcal{M}_{01} +\mathcal{H}_0\right) \check{\hat{X}}+\left(\mathcal{N}_{00} +\mathcal{B}_{00}\right)\check{Y}+\mathcal{C}_{00} \check{\hat{Y}}+\left(\mathcal{N}_{01}+ \mathcal{B}_{01}\right) \check{Z}_1+\mathcal{C}_{01} \check{\hat{Z}}_1\right. \\
&\quad \left.+\left(\mathcal{N}_{02} +\mathcal{B}_{02}\right)\check{Z}_2 +\mathcal{C}_{02}\check{\hat{Z}}_2+\left(\mathcal{N}_{03}+\mathcal{B}_{03}\right) \check{Z}_3+\mathcal{C}_{03} \check{\hat{Z}}_3\right]dt \\
&\quad +\sum_{i=2,3}\left[\mathcal{A}_i \check{X}+\left(\mathcal{H}_i +\mathcal{M}_{i 1}\right)\check{\hat{X}}+\left(\mathcal{N}_{i0} +\mathcal{B}_{i 0}\right)\check{Y}+\mathcal{C}_{i 0} \check{\hat{Y}}+\left(\mathcal{N}_{i1}+\mathcal{B}_{i1}\right) \check{Z}_1\right. \\
&\quad \left. +\mathcal{C}_{i1} \check{\hat{Z}}_1+\left(\mathcal{N}_{i 2}+\mathcal{B}_{i 2}\right) \check{Z}_2+\mathcal{C}_{i 2} \check{\hat{Z}}_2+\left(\mathcal{N}_{i 3} +\mathcal{B}_{i 3}\right) \check{Z}_3
+\mathcal{C}_{i 3} \check{\hat{Z}}_3 \right]dW_i, \quad t\in[0,T],\\
\check{X}(0) &=\mathcal{X}_0,
\end{aligned}\right.
\end{equation}
and
\begin{equation}\label{hat check X}
\left\{\begin{aligned}
d \check{\hat{X}}(t)&=\left[\left(\mathcal{A}_0+\mathcal{M}_{01} +\mathcal{H}_0\right) \check{\hat{X}}+\left(\mathcal{N}_{00} +\mathcal{B}_{00}+\mathcal{C}_{00}\right) \check{\hat{Y}}+\left(\mathcal{N}_{01} + \mathcal{B}_{01}+\mathcal{C}_{01}\right) \check{\hat{Z}}_1\right. \\
&\quad \left.+\left(\mathcal{N}_{02} +\mathcal{B}_{02}+\mathcal{C}_{02}\right) \check{\hat{Z}}_2+\left(\mathcal{N}_{03}+\mathcal{B}_{03}+\mathcal{C}_{03}\right) \check{\hat{Z}}_3\right]dt \\
&\quad +\left[\left(\mathcal{A}_3+\mathcal{M}_{3 1}+\mathcal{H}_3 \right) \check{\hat{X}}+\left(\mathcal{N}_{30}+\mathcal{B}_{3 0}+\mathcal{C}_{3 0}\right) \check{\hat{Y}}+\left(\mathcal{N}_{31}+ \mathcal{B}_{31}+\mathcal{C}_{31}\right) \check{\hat{Z}}_1\right. \\
&\quad \left. +\left(\mathcal{N}_{3 2} +\mathcal{B}_{3 2}+\mathcal{C}_{3 2}\right) \check{\hat{Z}}_2+\left(\mathcal{N}_{3 3}+\mathcal{B}_{3 3}+\mathcal{C}_{3 3}\right) \check{\hat{Z}}_3 \right]dW_3, \quad t\in[0,T],\\
\check{\hat{X}}(0) &=\mathcal{X}_0,
\end{aligned}\right.
\end{equation}
respectively. Observing the terminal condition of (\ref{Optimality system of leaders}), we put
\begin{equation}\label{Relation of Y and X}
Y(t)=P_1(t)X(t)+P_2(t) \hat{X}(t)+P_3(t) \check{X}(t)+P_4(t) \check{\hat{X}}(t),\quad t\in[0,T],
\end{equation}
where $P_i(\cdot) \in \mathbb{R}^{4 \times 5}, i=1,2,3,4$ are deterministic, differentiable unctions, satisfying $P_1(T)=\mathcal{G},P_i=O_{4 \times 5},i=2,3,4$.
Applying $\mathbb{E}\left[\cdot \mid \mathcal{G}_t^1\right]$, $\mathbb{E}\left[\cdot \mid \mathcal{G}_t^2\right]$ and $\mathbb{E}\left[\mathbb{E} \left[ \cdot \mid \mathcal{G}_t^1\right] \mid \mathcal{G}_t^2 \right]$ to (\ref{Relation of Y and X}), we obtain
\begin{equation}\label{hat Y}
\hat{Y}=\left(P_1 +P_2\right) \hat{X}+\left(P_3+P_4\right) \check{\hat{X}},
\end{equation}
\begin{equation}\label{check Y}
\check{Y}=\left(P_1 +P_3 \right) \check{X} +\left(P_2+P_4\right) \check{\hat{X}},
\end{equation}
\begin{equation}\label{hat and check Y}
\check{\hat{Y}}=\left(P_1+P_2 +P_3 +P_4\right) \check{\hat{X}},
\end{equation}
respectively. Applying It\^{o}'s formula to (\ref{Relation of Y and X}), we get
$$
\begin{aligned}
dY& =\bigg\{\big(\dot{P}_1+P_1 \mathcal{A}_0+P_1^2 \mathcal{N}_{00}\big) X+\Big[\dot{P}_2+P_1 \mathcal{M}_{01}+P_1 \mathcal{N}_{00} P_2+P_1 \mathcal{B}_{00}\left(P_1+P_3\right) \\
&\quad +P_2\left(\mathcal{A}_0+\mathcal{M}_{01}\right)+P_2 \mathcal{N}_{00}\left(P_1+P_2\right)\Big] \hat{X}+\Big[\dot{P}_3+P_1 \mathcal{N}_{00} P_3+P_1 \mathcal{B}_{00}\left(P_1+P_3\right)+P_3 \mathcal{A}_0 \\
&\quad +P_3\left(\mathcal{N}_{00}+\mathcal{B}_{00}\right)\left(P_1+P_3\right)\Big] \check{X}+\Big[\dot{P}_4+P_1 \mathcal{H}_0+P_1 \mathcal{N}_{00} P_4+P_1 \mathcal{B}_{00}\left(P_2+P_4\right) \\
&\quad +P_1 \mathcal{C}_{00}\left(P_1+P_2+P_3+P_4\right)+P_2 \mathcal{H}_0+P_2 \mathcal{N}_{00}\left(P_3+P_4\right)+P_2\left(\mathcal{B}_{00}+\mathcal{C}_{00}\right)\left(P_1+P_2 \right. \\
&\quad \left.+P_3+P_4\right)+P_3\left(\mathcal{M}_{01}+\mathcal{H}_0\right)+P_3\left(\mathcal{N}_{00}+\mathcal{B}_{00}\right)\left(P_2+P_4\right)+P_3 \mathcal{C}_{00}\left(P_1+P_2+P_3+P_4\right) \\
&\quad +P_4\left(\mathcal{A}_0+\mathcal{M}_{01}+\mathcal{H}_0\right)+P_4\left(\mathcal{N}_{00}+\mathcal{B}_{00}+\mathcal{C}_{00}\right)\left(P_1+P_2+P_3+P_4\right)\Big] \check{\hat{X}} \\
&\quad +\sum\limits_{i=1}^3 \left[ P_1 \mathcal{N}_{0i} Z_i+P_2 \mathcal{N} _{0i}\hat{Z}_i+\big(P_1 \mathcal{B}_{0i}+P_3\mathcal{N}_{0i}+P_3\mathcal{B}_{0i}\big)\check{Z}_i \right.\\
&\quad \left.+\big[P_1 \mathcal{C}_{0i}+P_2\left(\mathcal{B}_{0i}+\mathcal{C}_{0i}\right)+P_3 C_{0i}+P_4\left(\mathcal{N}_{0i}+\mathcal{B}_{0i}+\mathcal{C}_{0i}\right)\big] \check{\hat{Z}}_i \right] \bigg\}dt  \\
&\quad +\bigg\{\left(P_1 \mathcal{A}_1+P_1^2 \mathcal{N}_{10}\right) X+\big[P_1 \mathcal{M}_{11}+P_1 \mathcal{N}_{10} P_2+P_2\left(\mathcal{A}_1+\mathcal{M}_{11}\right)
 +P_2 \mathcal{N}_{10}\left(P_1+P_2\right)\big] \hat{X} \\
&\quad +\big[P_1 \mathcal{N}_{10} P_3+P_1 \mathcal{B}_{10}\left(P_1+P_3\right)\big] \check{X}+\big[P_1 \mathcal{H}_1+P_1 \mathcal{N}_{10} P_4+P_1 \mathcal{B}_{10}\left(P_2+P_4\right) \\
&\quad +P_1 \mathcal{C}_{10}\left(P_1+P_2+P_3+P_4\right)+P_2 \mathcal{H}_1+P_2 \mathcal{N}_{10}\left(P_3+P_4\right)+P_2\left(\mathcal{B}_{10}+\mathcal{C}_{10}\right)\left(P_1+P_2\right. \\
&\quad \left.+P_3 +P_4\right)\big] \check{\hat{X}}+\sum\limits_{i=1}^3 \left[P_1 \mathcal{N}_{1i} Z_i+P_2 \mathcal{N}_{1i} \hat{Z}_i+P_1 \mathcal{B}_{1i} \check{Z}_i+\left(P_1 \mathcal{C}_{1i}+P_2 \mathcal{B}_{1i}+P_2 \mathcal{C}_{1i}\right) \check{\hat{Z}}_i\right] \bigg\}dW_1\\
&\quad + \bigg\{\left(P_1 \mathcal{A}_2+P_1^2 \mathcal{N}_{20}\right) X+\left(P_1 \mathcal{M}_{21}+P_1 \mathcal{N}_{20} P_2\right) \hat{X}+\big[P_1 \mathcal{N}_{20} P_3+P_1 \mathcal{B}_{20}\left(P_1+P_3\right)
+P_3 \mathcal{A}_2 \\
&\quad +P_3\left(\mathcal{N}_{20}+\mathcal{B}_{20}\right)\left(P_1+P_3\right)\big] \check{X}+\big[P_1 \mathcal{H}_2+P_1 \mathcal{N}_{20} P_4+P_1 \mathcal{B}_{20}\left(P_2+P_4\right)+P_1 \mathcal{C}_{20}\left(P_1+P_2 \right. \\
&\quad \left. +P_3+P_4\right)+P_3\left(\mathcal{M}_2+\mathcal{H}_2\right)+P_3\left(\mathcal{N}_{20}+\mathcal{B}_{20}\right)\left(P_2+P_4\right)+P_3 \mathcal{C}_{20}\left(P_1+P_2+P_3+P_4\right)\big] \check{\hat{X}} \\
&\quad +\sum\limits_{i=1}^3 \Big[ P_1 \mathcal{N}_{2i} Z_i+\left(P_1 \mathcal{B}_{2i}+P_3 \mathcal{N}_{2i}+P_3 \mathcal{B}_{2i}\right) \check{Z}_i
 +\left(P_1 \mathcal{C}_{2i}+P_3 \mathcal{C}_{2i}\right)  \check{\hat{Z}}_t^i\Big]\bigg\} dW_2 \\
&\quad +\bigg\{\left(P_1 \mathcal{A}_3+P_1^2 \mathcal{N}_{30}\right) X+\big[P_1 \mathcal{M}_{31}+P_1 \mathcal{N}_{30}P_2+P_2\left(\mathcal{A}_3+\mathcal{M}_{31}\right)+P_2 \mathcal{N}_{30}\left(P_1+P_2\right)\big]\hat{X} \\
&\quad +\big[P_1 \mathcal{N}_{30} P_3+P_1 \mathcal{B}_{30}\left(P_1+P_3\right)+P_3 \mathcal{A}_3+P_3\left(\mathcal{N}_{30}+\mathcal{B}_{30}\right)\left(P_1+P_3\right)\big] \check{X} +\big[P_1 \mathcal{H}_3 \\
&\quad \left. +P_1 \mathcal{N}_{30} P_4+P_1 \mathcal{B}_{30}\left(P_2+P_4\right)+P_1 \mathcal{C}_{30}\left(P_1+P_2+P_3+P_4\right)+P_2 \mathcal{H}_3+P_2 \mathcal{N}_{30}\left(P_3+P_4\right) \right. \\
&\quad +P_2\left(\mathcal{B}_{30}+\mathcal{C}_{30}\right)\left(P_1+P_2+P_3+P_4\right) +P_3\left(\mathcal{M}_{31}+\mathcal{H}_3\right)+P_3\left(\mathcal{N}_{30}+\mathcal{B}_{30}\right)\left(P_2+P_4\right) \\
&\quad +P_3 \mathcal{C}_{30}\left(P_1+P_2+P_3+P_4\right)+P_4\left(\mathcal{A}_3+\mathcal{M}_{31}+\mathcal{H}_3\right) +P_4\left(\mathcal{N}_{30}+\mathcal{B}_{30}+\mathcal{C}_{30}\right)  \\
&\quad \times\left(P_1+P_2+P_3+P_4\right)\big] \check{\hat{X}} +\sum\limits_{i=1}^3 \left[P_1 \mathcal{N}_{3i} Z_t^{i}+P_2 \mathcal{N}_{3i} \hat{Z}_t^{i}+\left(P_1 \mathcal{B}_{3i}  +P_3 \mathcal{N}_{3i}+P_3 \mathcal{B}_{3i}\right) \check{Z}_t^{i}\right.  \\
&\quad \left.+\big[P_1 \mathcal{C}_{3i}+P_2\left(\mathcal{B}_{3i}+\mathcal{C}_{3i}\right)+P_3 \mathcal{C}_{3i}+P_4\left(\mathcal{N}_{3i}+\mathcal{B}_{3i}+\mathcal{C}_{3i}\right)\big] \check{\hat{Z}}_i \right] \bigg\}dW_3\\
\end{aligned}
$$
\begin{equation}\label{Apply Ito formula leaders}
\begin{aligned}
& =-\bigg\{\left(\mathcal{Q}+\overline{\mathcal{A}}_0 P_1\right) X+\big[\overline{\mathcal{A}}_{0}P_2+\overline{\mathcal{M}}_{01}\left(P_1+P_2\right)\big]\hat{X}+\overline{\mathcal{A}}_{0}P_3 \check{X} +\big[\check{\hat{\mathcal{Q}}}+\overline{\mathcal{A}}_{0}P_4 \\
&\quad +\overline{\mathcal{M}}_{01}\left(P_3+P_4\right)+\mathcal{I}_0\left(P_1+P_2+P_3+P_4\right)\big] \check{\hat{X}} +\overline{\mathcal{A}}_1 Z_1 +\overline{\mathcal{M}}_{11} \hat{Z}_1 +\mathcal{I}_1 \check{\hat{Z}}_1+\overline{\mathcal{A}}_2 Z_2   \\
&\quad +\overline{\mathcal{M}}_{21} \hat{Z}_2+\mathcal{I}_2 \check{\hat{Z}}_2+\overline{A}_3 Z_3 +\overline{M}_{31} \hat{Z}_3+\mathcal{I}_3 \check{\hat{Z}}_3 \bigg\}dt +\sum_{i=1}^3 Z_i d W_i.
\end{aligned}
\end{equation}
Comparing the diffusion terms on both sides of (\ref{Apply Ito formula leaders}), we obtain
\begin{equation}\label{coupled Z1}
\begin{aligned}
Z_1 & =\left(P_1 \mathcal{A}_1+P_1^2 \mathcal{N}_{10}\right) X+\big[P_1 \mathcal{M}_{11}+P_1 \mathcal{N}_{10} P_2+P_2\left(\mathcal{A}_1+\mathcal{M}_{11}\right)+P_2 \mathcal{N}_{10}\left(P_1+P_2\right)\big]\hat{X}  \\
&\quad +\big[P_1 \mathcal{N}_{10} P_3+P_1 \mathcal{B}_{10}\left(P_1+P_3\right)\big] \check{X}+\big[P_1 \mathcal{H}_1+P_1 \mathcal{N}_{10} P_4+P_1 \mathcal{B}_{10}\left(P_2+P_4\right) \\
&\quad +P_1 \mathcal{C}_{10}\left(P_1+P_2+P_3+P_4\right)+P_2 \mathcal{H}_1+P_2 \mathcal{N}_{10}\left(P_3+P_4\right)+P_2\left(\mathcal{B}_{10}+\mathcal{C}_{10}\right)\left(P_1+P_2 \right. \\
&\quad \left.+P_3+P_4\right)\big] \check{\hat{X}}+\sum\limits_{i=1}^3 \left[P_1 \mathcal{N}_{1i} Z_i+P_2 \mathcal{N}_{1i} \hat{Z}_t^i+P_1 \mathcal{B}_{1i} \check{Z}_i+\left(P_1 \mathcal{C}_{1i}+P_2 \mathcal{B}_{1i}+P_2 \mathcal{C}_{1i}\right) \check{\hat{Z}}_i\right],
\end{aligned}
\end{equation}
\begin{equation}\label{coupled Z2}
\begin{aligned}
Z_2 & = \left(P_1 \mathcal{A}_2+P_1^2 \mathcal{N}_{20}\right) X_t+\left(P_1 \mathcal{M}_{21}+P_1 \mathcal{N}_{20} P_2\right) \hat{X}_t+\big[P_1 \mathcal{N}_{20} P_3+P_1 \mathcal{B}_{20}\left(P_1+P_3\right)+P_3 \mathcal{A}_2 \\
&\quad +P_3\left(\mathcal{N}_{20}+\mathcal{B}_{20}\right)\left(P_1+P_3\right)\big] \check{X} +\big[P_1 \mathcal{H}_2+P_1 \mathcal{N}_{20} P_4+P_1 \mathcal{B}_{20}\left(P_2+P_4\right)+P_1 \mathcal{C}_{20}\left(P_1+P_2\right. \\
&\quad \left. +P_3+P_4\right)+P_3\left(\mathcal{M}_2+\mathcal{H}_2\right)+P_3\left(\mathcal{N}_{20}+\mathcal{B}_{20}\right)\left(P_2+P_4\right)+P_3 \mathcal{C}_{20}\left(P_1+P_2+P_3+P_4\right)\big] \check{\hat{X}} \\
&\quad +\sum\limits_{i=1}^3 \Big[ P_1 \mathcal{N}_{2i} Z_i+\left(P_1 \mathcal{B}_{2i}+P_3 \mathcal{N}_{2i}+P_3 \mathcal{B}_{2i}\right) \check{Z}_i+\left(P_1 \mathcal{C}_{2i}+P_3 \mathcal{C}_{2i}\right)  \check{\hat{Z}}_i\Big],
\end{aligned}
\end{equation}
\begin{equation}\label{coupled Z3}
\begin{aligned}
Z_3 & =\left(P_1 \mathcal{A}_3+P_1^2 \mathcal{N}_{30}\right) X+\big[P_1 \mathcal{M}_{31}+P_1 \mathcal{N}_{30}P_2+P_2\left(\mathcal{A}_3+\mathcal{M}_{31}\right)+P_2 \mathcal{N}_{30}\left(P_1+P_2\right)\big]\hat{X} \\
&\quad +\big[P_1 \mathcal{N}_{30} P_3+P_1 \mathcal{B}_{30}\left(P_1+P_3\right)+P_3 \mathcal{A}_3+P_3\left(\mathcal{N}_{30}+\mathcal{B}_{30}\right)\left(P_1+P_3\right)\big] \check{X} +\big[P_1 \mathcal{H}_3 \\
&\quad +P_1 \mathcal{N}_{30} P_4+P_1 \mathcal{B}_{30}\left(P_2+P_4\right)+P_1 \mathcal{C}_{30}\left(P_1+P_2+P_3+P_4\right)+P_2 \mathcal{H}_3+P_2 \mathcal{N}_{30}\left(P_3+P_4\right) \\
&\quad +P_2\left(\mathcal{B}_{30}+\mathcal{C}_{30}\right)\left(P_1+P_2+P_3+P_4\right) +P_3\left(\mathcal{M}_{31}+\mathcal{H}_3\right)+P_3\left(\mathcal{N}_{30}+\mathcal{B}_{30}\right)\left(P_2+P_4\right) \\
&\quad +P_3 \mathcal{C}_{30}\left(P_1+P_2+P_3+P_4\right)+P_4\left(\mathcal{A}_3+\mathcal{M}_{31}+\mathcal{H}_3\right) +P_4\left(\mathcal{N}_{30}+\mathcal{B}_{30}+\mathcal{C}_{30}\right)  \\
&\quad \times\left(P_1+P_2+P_3+P_4\right)\big] \check{\hat{X}} +\sum\limits_{i=1}^3 \left[P_1 \mathcal{N}_{3i} Z_i+P_2 \mathcal{N}_{3i} \hat{Z}_i+\left(P_1 \mathcal{B}_{3i}  +P_3 \mathcal{N}_{3i}+P_3 \mathcal{B}_{3i}\right) \check{Z}_i\right.  \\
&\quad \left.+\big[P_1 \mathcal{C}_{3i}+P_2\left(\mathcal{B}_{3i}+\mathcal{C}_{3i}\right)+P_3 \mathcal{C}_{3i}+P_4\left(\mathcal{N}_{3i}+\mathcal{B}_{3i}+\mathcal{C}_{3i}\right)\big] \check{\hat{Z}}_i \right],
\end{aligned}
\end{equation}
respectively. However, since $Z_1(\cdot),Z_2(\cdot),Z_3(\cdot)$ are coupled in (\ref{coupled Z1}), (\ref{coupled Z2}) and (\ref{coupled Z3}), we have to add some additional assumptions to obtain the explicit expressions for them.

Thanks to Shi et al. \cite{Shi-Wang-Xiong-2020}, we could obtain the explicit expression $Z_1(\cdot),Z_2(\cdot),Z_3(\cdot)$:
\begin{equation}\label{Z--}
Z_i = \mathcal{N}_i X+\hat{\mathcal{N}}_i \hat{X} +\check{\mathcal{N}}_i \check{X}+\check{\hat{\mathcal{N}}}_i \check{\hat{X}}, \quad i=1,2,3,
\end{equation}
together with
\begin{equation}\label{hatcheck and check and hat Z}
\check{\hat{Z}}_i= N_i \check{\hat{X}}, \quad \hat{Z}_i =\hat{N}_i \hat{X}+\check{\hat{N}}_i \check{\hat{X}},\quad
\check{Z}_i =\check{N}_i \hat{X}+\tilde{N}_i \check{\hat{X}}, \quad i=1,2,3,
\end{equation}
by decoupling (\ref{coupled Z1}), (\ref{coupled Z2}) and (\ref{coupled Z3}) with some filtering technique. The detail is given in the Appendix, together with the definitions of $N_i,\hat{N}_i,\check{N}_i,\tilde{N}_i,\check{\hat{N}}_i,\mathcal{N}_i,\hat{\mathcal{N}}_i,\check{\mathcal{N}}_i,\check{\hat{\mathcal{N}}}_i$, $i=1,2,3$, and some assumptions.

Next, comparing the drift terms in (\ref{Apply Ito formula leaders}) and substituting (\ref{Z--}), (\ref{hatcheck and check and hat Z}) into it, we obtain
\begin{equation}\label{Riccation equations of leaders P1}
\begin{aligned}
&\dot{P}_1+P_1 \mathcal{A}_0+\overline{\mathcal{A}}_0 P_1+P_1^2 \mathcal{N}_{00}+\mathcal{Q}+\sum_{i=1}^3\left(P_1 \mathcal{N}_{0 i}+\overline{\mathcal{A}}_i\right) \mathcal{N}_i=O_{4 \times 5},\quad P_1(T)=\mathcal{G},
\end{aligned}
\end{equation}
\begin{equation}\label{Riccation equations of leaders P2}
\begin{aligned}
&\dot{P}_2+P_1\mathcal{M}_{01}+2P_1 \mathcal{N}_{00} P_2+P_1 \mathcal{B}_{00}\left(P_1+P_3\right)+P_2\left(\mathcal{A}_0+\mathcal{M}_{01}\right)+P_2^2 \mathcal{N}_{00}+\overline{\mathcal{A}}_0 P_2+\overline{\mathcal{M}}_{01}\left(P_1\right.\\
& \left.+P_2\right)+\sum_{i=1}^3\Big[\left(P_1 \mathcal{N}_{0 i}+\overline{A}_i\right) \hat{\mathcal{N}}_i
+\left(P_2 \mathcal{N}_{01}+\overline{\mathcal{M}}_{i 1}\right)\big(\mathcal{N}_i+\hat{\mathcal{N}}_i\big)\Big] =O_{4 \times 5},\quad P_2(T)=O_{4 \times 5},
\end{aligned}
\end{equation}
\begin{equation}\label{Riccation equations of leaders P3}
\begin{aligned}
& \dot{P}_3+P_1 \mathcal{N}_{00} P_3+P_1 \mathcal{B}_{00}\left(P_1+P_3\right)+P_3 \mathcal{A}_0+P_3\left(\mathcal{N}_{00}+\mathcal{B}_{00}\right)\left(P_1+P_3\right)+\overline{\mathcal{A}}_0 P_3 \\
& +\sum_{i=1}^3\Big[\left(P_1 \mathcal{N}_{0 i}+\overline{\mathcal{A}}_i\right) \check{N}_i+\left(P_1 \mathcal{B}_{0 i}+P_3 \mathcal{N}_{0 i}+P_3 \mathcal{B}_{0 i}\right)\left(\mathcal{N}_i+\check{\mathcal{N}}_i\right)\Big]
 =O_{4 \times 5},\quad P_3(T)=O_{4 \times 5},
\end{aligned}
\end{equation}
\begin{equation}\label{Riccation equations of leaders P4}
\begin{aligned}
& \dot{P}_4+\left(P_1+P_2\right) \mathcal{H}_0+P_1 \mathcal{N}_{00} P_4+P_1 \mathcal{B}_{30}\left(P_2+P_4\right)+P_2 \mathcal{N}_{00}\left(P_3+P_4\right) \\
& +P_3\left(\mathcal{M}_0+\mathcal{H}_0\right)+P_3\left(\mathcal{N}_{00}+\mathcal{B}_{00}\right)\left(P_2+P_4\right)+P_4\left(\mathcal{A}_0+\mathcal{M}_{01}+\mathcal{H}_0\right)\\
& +\big[\mathcal{I}_0+\left(P_1+P_3\right) \mathcal{C}_{00}+P_2\left(\mathcal{B}_{00}+\mathcal{C}_{00}\right)+P_4\left(\mathcal{N}_{00}+\mathcal{B}_{00}+\mathcal{C}_{00}\right)\big]\left(P_1+P_2+P_3+P_4\right) \\
& +\check{\hat{\mathcal{Q}}}+\overline{\mathcal{A}}_0 P_4+\overline{\mathcal{M}}_{01}\left(P_3+P_4\right)+\sum_{i=1}^3\left[\left(P_1 \mathcal{N}_{0 i}+\overline{A}_i\right) \check{\hat{N}}_i+\left(P_2 \mathcal{N}_{0 i}+\mathcal{M}_{i 1}\right)\big(\check{N}_i+\check{\hat{N}}_i\big)\right. \\
& +\left(P_1 \mathcal{B}_{0 i}+P_3 \mathcal{N}_{0 i}+P_3 \mathcal{B}_{0 i}\right)\big(\hat{\mathcal{N}}_i+\check{\hat{\mathcal{N}}}_i\big)+\big[P_1 \mathcal{C}_{0 i}+P_2\left(\mathcal{B}_{0 i}+\mathcal{C}_{0 i}\right)+P_3 \mathcal{C}_{0 i} \\
& \left.+P_4\left(\mathcal{N}_{0 i}+\mathcal{B}_{0 i}+\mathcal{C}_{0 i}\right)+\mathcal{I}_i\big]\big(\mathcal{N}_i+\hat{\mathcal{N}}_i+\mathcal{N}_i+\check{\hat{\mathcal{N}}}_i\big)\right] =O_{4 \times 5},\quad P_4(T)=O_{4 \times 5}.
\end{aligned}
\end{equation}

Since $\mathcal{N}_i,\hat{\mathcal{N}}_i,\check{\mathcal{N}}_i,\check{\hat{\mathcal{N}}}_i$ depend on $P_i,i=1,2,3,4$, the solvability of this coupled system of asymmetric Riccati equations is rather difficult. We will discuss this problem in the future.

Finally, substituting (\ref{check Y}), (\ref{hat and check Y}), (\ref{hatcheck and check and hat Z}) into (\ref{SMP of leaders-}), we obtain the state-estimate feedback form of the Nash equilibrium point of the leaders' problem as follows:
\begin{equation*}
\begin{aligned}
u_3^* & =-R_3^{-1}\left\{\left[\overline{\mathcal{D}}_0^\top\left(P_1+P_3\right)+\sum_{i=1}^3 \overline{\mathcal{D}}_i^\top\left(\mathcal{N}_i+\check{\mathcal{N}}_i\right)\right] \check{X}\right. \\
&\quad +\left[\mathcal{F}_1^\top+\overline{\mathcal{D}}_0^\top\left(P_2+P_4\right)+\overline{\mathcal{M}}_{04}^\top\left(P_1+P_2+P_3+P_4\right)+\sum_{i=1}^3 \overline{\mathcal{D}}_i^\top\big(\hat{\mathcal{N}}_i+\check{\hat{\mathcal{N}}}_i\big)\right. \\
&\quad \left.\left.+\sum_{i=1}^3 \overline{\mathcal{M}}_{i 4}^\top\left(\mathcal{N}_i+\hat{\mathcal{N}}_i+\check{\mathcal{N}}_i+\check{\hat{\mathcal{N}}}_i\right)\right] \check{\hat{X}}\right\},  \\
u_4^* & =-R_4^{-1}\left\{\left[\overline{\mathcal{E}}_0^{\top}\left(P_1+P_3\right)+\sum_{i=1}^3 \overline{\mathcal{E}}_i^{\top}\left(\mathcal{N}_i+\check{\mathcal{N}}_i\right)\right] \check{X}\right. \\
\end{aligned}
\end{equation*}
\begin{equation}\label{Nash equilibrium point of the leaders}
\begin{aligned}
&\quad +\left[\mathcal{F}_2^\top+\overline{\mathcal{E}}_0^\top\left(P_2+P_4\right)+\overline{\mathcal{M}}_{05}^\top\left(P_1+P_2+P_3+P_4\right)+\sum_{i=1}^3 \overline{\mathcal{E}}_i^\top\big(\hat{\mathcal{N}}_i+\check{\hat{\mathcal{N}}}_i\big)\right. \\
&\quad \left.\left.+\sum_{i=1}^3 \overline{\mathcal{M}}_{i 5}^\top\big(\mathcal{N}_i+\hat{\mathcal{N}}_i+\check{\mathcal{N}}_i+\check{\hat{\mathcal{N}}}_i\big)\right] \check{\hat{X}}\right\}.
\end{aligned}
\end{equation}

We summarize the above argument in the following theorem.
\begin{thm}\label{thm4}
Let Assumptions \ref{A1}, \ref{A2}, \ref{A4}-\ref{A8} hold, $P_i(\cdot)$, $i=1,2,3,4$ satisfy the system (\ref{Riccation equations of leaders P1})-(\ref{Riccation equations of leaders P4}) of Riccati equations. Let $(\hat{X}(\cdot), \check{X}(\cdot), \check{\hat{X}}(\cdot))$ be the adapted solutions to (\ref{hat X}), (\ref{check X}), (\ref{hat check X}), respectively, with relations (\ref{hat Y}), (\ref{check Y}), (\ref{hat and check Y}), (\ref{hatcheck and check and hat Z}). Then $\left(u_3^*(\cdot),u_4^*(\cdot)\right)$ given by (\ref{Nash equilibrium point of the leaders}) is the state-estimate feedback representation for the Nash equilibrium point of the leaders.
\end{thm}

Using similar method, we can derive a ``non-anticipating" form about the Nash equilibrium point of the followers. In fact, by (\ref{Nash equilibrium point of the followers}), (\ref{X,Y,Z}), (\ref{hat Y}), (\ref{hat and check Y}), (\ref{hatcheck and check and hat Z}), we obtain
\begin{equation*}
\begin{aligned}
u_1^* & =\bigg\{\left(\begin{array}{ccc}
			-L_{11} & 0 & 0
		\end{array}\right)+\left(\begin{array}{ccc}
		\alpha_1 & 0 & 0
		\end{array}\right)\left(P_1+P_2\right) +\left(\begin{array}{ccc}
		\beta_1 & 0 & 0
		\end{array}\right)\left(\mathcal{N}_1+\hat{\mathcal{N}}_1\right)\\
&\quad +\left(\begin{array}{ccc}
		\gamma_1 & 0 & 0
		\end{array}\right)\left(\mathcal{N}_3+\hat{\mathcal{N}}_3\right)\bigg\}\hat{X}+\bigg\{L_{12}R_3^{-1}\bigg[\overline{\mathcal{D}}_0^\top\left(P_1+P_3\right)
 +\sum_{i=1}^3\overline{\mathcal{D}}_i^{\top}\left(\mathcal{N}_i+\check{\mathcal{N}}_i\right)  \\
&\quad +\mathcal{F}_1^\top+\overline{\mathcal{D}}_0^\top\left(P_2+P_4\right)+\overline{\mathcal{M}}_{04}^\top\left(P_1+P_2+P_3+P_4\right)
 +\sum_{i=1}^3 \overline{\mathcal{D}}_i^{\top}\big(\hat{\mathcal{N}}_i+\check{\hat{\mathcal{N}}}_i\big) \\
&\quad +\sum_{i=1}^3 \overline{\mathcal{M}}_{i 4}^\top\big(\mathcal{N}_i+\hat{\mathcal{N}}_i+\check{\mathcal{N}}_i
 +\check{\hat{\mathcal{N}}}_i\big)\bigg] +L_{13}R_4^{-1}\bigg[\overline{\mathcal{E}}_0^\top\left(P_1+P_3\right)+\sum_{i=1}^3 \overline{\mathcal{E}}_i^\top\big(\mathcal{N}_i+\check{\mathcal{N}}_i\big) \\
&\quad +\mathcal{F}_2^{\top}+\overline{\mathcal{E}}_0^\top\left(P_2+P_4\right)+\overline{\mathcal{M}}_{05}^\top\left(P_1+P_2+P_3+P_4\right)
 +\sum_{i=1}^3 \overline{\mathcal{E}}_i^{\top}\big(\hat{\mathcal{N}}_i+\check{\hat{\mathcal{N}}}_i\big)\\
&\quad +\sum_{i=1}^3 \overline{\mathcal{M}}_{i 5}^\top\big(\mathcal{N}_i+\hat{\mathcal{N}}_i+\check{\mathcal{N}}_i+\check{\hat{\mathcal{N}}}_i\big)\bigg] +\left(\begin{array}{ccc}
			\alpha_1 & 0 & 0
		\end{array}\right)\left(P_3+P_4\right)\\
&\quad +\left(\begin{array}{ccc}
			\beta_1 & 0 & 0
		\end{array}\right)\big(\check{\mathcal{N}}_1+\check{\hat{\mathcal{N}}}_1\big)+\left(\begin{array}{ccc}
			\gamma_1 & 0 & 0
		\end{array}\right)\big(\check{\mathcal{N}}_3+\check{\hat{\mathcal{N}}}_3\big)\bigg\}\check{\hat{X}}, \\
u_2^* & =\bigg\{\left(\begin{array}{ccc}
			-L_{21} & 0 & 0
		\end{array}\right)+\left(\begin{array}{ccc}
			\alpha_2 & 0 & 0
		\end{array}\right)\left(P_1+P_2\right) +\left(\begin{array}{ccc}
			\beta_2 & 0 & 0
		\end{array}\right)\big(\mathcal{N}_1+\hat{\mathcal{N}}_1\big)\\
&\quad +\left(\begin{array}{ccc}
			\gamma_2 & 0 & 0
		\end{array}\right)\big(\mathcal{N}_3+\hat{\mathcal{N}}_3\big)\bigg\}\hat{X}+\bigg\{L_{22}R_3^{-1}\bigg[\overline{\mathcal{D}}_0^\top\left(P_1+P_3\right)
 +\sum_{i=1}^3\overline{\mathcal{D}}_i^\top\left(\mathcal{N}_i+\check{\mathcal{N}}_i\right) \\
&\quad +\mathcal{F}_1^\top+\overline{\mathcal{D}}_0^{\top}\left(P_2+P_4\right)+\overline{\mathcal{M}}_{04}^\top\left(P_1+P_2+P_3+P_4\right)
 +\sum_{i=1}^3 \overline{\mathcal{D}}_i^\top\big(\hat{\mathcal{N}}_i+\check{\hat{\mathcal{N}}}_i\big) \\
&\quad +\sum_{i=1}^3 \overline{\mathcal{M}}_{i 4}^\top\big(\mathcal{N}_i+\hat{\mathcal{N}}_i
 +\check{\mathcal{N}}_i+\check{\hat{\mathcal{N}}}_i\big)\bigg] +L_{23}R_4^{-1}\bigg[\overline{\mathcal{E}}_0^\top\left(P_1+P_3\right)
 +\sum_{i=1}^3 \overline{\mathcal{E}}_i^\top\big(\mathcal{N}_i+\check{\mathcal{N}}_i\big) \\
\end{aligned}
\end{equation*}
\begin{equation}\label{non-anticipating follower}
\begin{aligned}
&\quad +\mathcal{F}_2^\top+\overline{\mathcal{E}}_0^\top\left(P_2+P_4\right)+\overline{\mathcal{M}}_{05}^\top\left(P_1+P_2+P_3+P_4\right)
 +\sum_{i=1}^3 \overline{\mathcal{E}}_i^\top\big(\hat{\mathcal{N}}_i+\check{\hat{\mathcal{N}}}_i\big) \\
&\quad +\sum_{i=1}^3 \overline{\mathcal{M}}_{i 5}^\top\big(\mathcal{N}_i+\hat{\mathcal{N}}_i+\check{\mathcal{N}}_i+\check{\hat{\mathcal{N}}}_i\big)\bigg] +\left(\begin{array}{ccc}
			\alpha_2 & 0 & 0
		\end{array}\right)\left(P_3+P_4\right)\\
&\quad +\left(\begin{array}{ccc}
			\beta_2 & 0 & 0
		\end{array}\right)\left(\check{\mathcal{N}}_1+\check{\hat{\mathcal{N}}}_1\right)+\left(\begin{array}{ccc}
			\gamma_2 & 0 & 0
		\end{array}\right)\left(\check{\mathcal{N}}_3+\check{\hat{\mathcal{N}}}_3\right)\bigg\}\check{\hat{X}},
\end{aligned}
\end{equation}
where
$$
\begin{aligned}
\alpha_1 &:=\frac{\left(\begin{array}{cc}
	\overline{R}_2 b_0 &  \widetilde{P}_1 c_0\left(\sum\limits_{i=1}^{3}b_i c_i \right)
\end{array}\right)}{|N|},\quad \alpha_2:=\frac{\left(\begin{array}{cc}
	\widetilde{P}_2 b_0\left(\sum\limits_{i=1}^{3}b_i c_i \right)    &   \overline{R}_1 c_0
\end{array}\right)}{|N|},\\
\beta_1 &:=\frac{\left(\begin{array}{cc}
		\overline{R}_2 b_1 & \overline{R}_2 b_3\end{array}\right)}{|N|},\quad \beta_2 =\frac{\left(\begin{array}{cc}
		\widetilde{P}_2 b_1\left(\sum\limits_{i=1}^{3}b_i c_i \right)    &   \widetilde{P}_2 b_3\left(\sum\limits_{i=1}^{3}b_i c_i \right)
	\end{array}\right)}{|N|},
\end{aligned}
$$
$$
\gamma_1:=\frac{\left(\begin{array}{cc}
		\widetilde{P}_1 c_1\left(\sum\limits_{i=1}^{3}b_i c_i \right)    &   \widetilde{P}_1 c_3\left(\sum\limits_{i=1}^{3}b_i c_i \right)\end{array}\right)}{|N|},\quad
\gamma_2:=\frac{\left(\begin{array}{cc}
		\overline{R}_1 c_1 & \overline{R}_1 c_3\end{array}\right)}{|N|}.
$$

Up to now, the Stackelberg-Nash equilibrium $(u_1^*(\cdot),u_2^*(\cdot),u_3^*(\cdot),u_4^*(\cdot))$ of the game is obtained, which is represented as the state-estimate feedback form in (\ref{non-anticipating follower}) and (\ref{Nash equilibrium point of the leaders}).

\section{Conclusion}

In this paper, inspired by the system of collective leadership and asymmetric information among players, we have discussed overlapping information LQ Stackelberg stochastic differential game with two leaders and two followers. There are three distinct features in our paper. Firstly, the diffusion terms in the state equation of the system contain both state variables and control variables. Secondly, there are different/asymmetric noisy information between the leaders and the followers, and they both need to solve partial information LQ nonzero-sum Nash games. Finally, new systems of coupled Riccati equation are introduced to give the state-estimate feedback representation of the Stackelberg-Nash equilibrium.

There are still some topics that need further research. The general solvability of the system of coupled Riccati equations (\ref{Riccation equations of leaders P1})-(\ref{Riccation equations of leaders P4}) is rather challenging. Potential applications in practice are interesting. Extensions of our problem to LQ Stackelberg stochastic differential games of mean-field type, mixed leadership model, time delayed model, and more complex information structure, are all our future researching interest.

%%%%%%%%%%%%%%%%%%%%%%%%%%%%%%%%%%%%%%%%%%%%%%%%%%%%%%%%%%%%

\section*{Appendix}

{\it Proof of Theorem \ref{thm1}.}\begin{proof}
{\it Necessity.} We define the followers' Hamiltonian functions as
\begin{equation*}\label{Hamilton of followers}
\begin{aligned}
& H_j\left(t, x, u_1, u_2, u_3, u_4, p_j, k_{j1},k_{j2},k_{j3}\right):= p_j\left(a_0 x+b_0 u_1+c_0 u_2+d_0 u_3+e_0 u_4\right)  \\
& + \sum_{i=1}^3 k_{j i}\left(a_i x+b_i u_1+c_i u_2+d_i u_3+e_i u_{4}\right) -\frac{1}{2} Q_j x^2-\frac{1}{2} R_j u_j^2, \quad j=1,2 .
\end{aligned}
\end{equation*}
Using the SMP with partial information (\cite{Baghery-Oksendal-2007} or \cite{Nie-Wang-Yu-2022}), we can obtain (\ref{SMP of followers}), (\ref{adjoint BSDEs of followers}).

{\it Sufficiency.} For any $u_1(\cdot) \in \mathcal{U}_1$, letting $x_1(\cdot):=x^{u_1, u_2^*, u_3, u_4}(\cdot)-x^{u_1^*, u_2^*, u_3, u_4}(\cdot)$, $x^{u_3, u_4}(\cdot):=x^{u_1^*, u_2^*, u_3, u_4}(\cdot)$, we have
$$
\begin{aligned}
& J_1\left(u_1(\cdot), u_2^*(\cdot), u_3(\cdot), u_4(\cdot)\right)-J_1\left(u_1^*(\cdot), u_2^*(\cdot), u_3(\cdot), u_4(\cdot)\right) \\
& =\frac{1}{2} \mathbb{E}\left[\int_0^T \left[Q_1 x_1^2+ R_1\left(u_1-u_1^*\right)^2\right] d t+ G_1 x_1(T)\right]+\Theta_f,
\end{aligned}
$$
where
$$
\Theta_f:=\mathbb{E}\left[\int_0^T \Big[Q_1 x_1 x^{u_3,u_4}+R_1 u_1^* \left(u_1-u_1^*\right)\Big] d t+ G_1  x_1(T) x^{u_3,u_4}(T)\right].
$$

Since
$$
\mathbb{E}\left[ G_1 x^{u_3,u_4}(T) x_1(T)\right]=-\mathbb{E}\left[p_1(T) x_1(T) \right],\quad x_1(0)=0,
$$
applying It\^{o}'s formula to $ x_1(\cdot) p_1(\cdot)$, we get
$$
\begin{aligned}
\mathbb{E}\left[x_1(T) p_1(T)\right] = \mathbb{E} \int_0^T\bigg[ b_0 \left(u_1-u_1^*\right) p_1+ Q_1 x_1 x^{u_3,u_4}+\sum\limits_{i=1}^3 b_i k_{1i}(u_1-u_1^*) \bigg] d t.
\end{aligned}
$$
By (\ref{adjoint BSDEs of followers}) and the properties of conditional expectation, we obtain
$$
\begin{aligned}
\Theta_f & =\mathbb{E}\left[\int_0^T \Big[Q_1 x_1 x^{u_3,u_4}+R_1 u_1^* \left(u_1-u_1^*\right)\Big] d t - x_1(T) p_1(T) \right] \\
& =\mathbb{E} \int_0^T\bigg[ R_1 u_1^*\left(u_1-u_1^*\right)- b_0 p_1\left(u_1-u_1^*\right)- \sum\limits_{i=1}^3 b_i k_{1i}(u_1-u_1^*)\bigg] d t = 0.
\end{aligned}
$$
Then by Assumption \ref{A1}, we get
$$
J_1\left(u_1(\cdot), u_2^*(\cdot), u_3(\cdot), u_4(\cdot)\right)-J_1\left(u_1^*(\cdot), u_2^*(\cdot), u_3(\cdot), u_4(\cdot)\right) \geq 0 .
$$

Similarly, for any $u_2(\cdot) \in \mathcal{U}_2$, we also have
$$
J_2\left(u_1^*(\cdot), u_2(\cdot), u_3(\cdot), u_4(\cdot)\right)-J_2\left(u_1^*(\cdot), u_2^*(\cdot), u_3(\cdot), u_4(\cdot)\right) \geq 0 .
$$
Therefore, $\left(u_1^*(\cdot), u_2^*(\cdot)\right)$ is a Nash equilibrium point for the followers' problem.
\end{proof}
{\it Proof of Lemma \ref{lemma1}.}\begin{proof}
Let $\tilde{P}(\cdot)=\tilde{P}_1(\cdot)+\tilde{P}_2(\cdot)$. It follows from Assumption \ref{A3} and (\ref{coupled Riccati equations of followers}) that
\begin{equation}\label{added Riccati equations of followers}
\dot{\tilde{P}}+2a_0 \tilde{P}+\sum\limits_{i=1}^3 a_i^2 \tilde{P}+Q_1+Q_2+\frac{\bar{B}_1 \tilde{P}^2}{R_1+\left(\sum\limits_{i=1}^3 b_i^2\right)\tilde{P}}=0, \quad
\tilde{P}(T)=G_1+G_2 .
\end{equation}
Obviously, (\ref{added Riccati equations of followers}) is a standard Riccati equation, it admits a unique solution $\tilde{P}(\cdot)$ (Chapter 6, \cite{Yong-Zhou-1999}). Introduce two auxiliary ordinary differential equations (ODEs):
\begin{equation}\label{auxiliary ODEs 1 of followers}
\dot{\check{\tilde{P}}}_1+2a_0 \check{\tilde{P}}_1+\sum\limits_{i=1}^3 a_i^2 \check{\tilde{P}}_1+Q_1+\frac{\bar{B}_1 \check{\tilde{P}}_1 \tilde{P}}{R_1+\left(\sum\limits_{i=1}^3 b_i^2\right)\tilde{P}}=0, \quad
\check{\tilde{P}}_1(T)=G_1,
\end{equation}
\begin{equation}\label{auxiliary ODEs 2 of followers}
\dot{\check{\tilde{P}}}_2+2a_0 \check{\tilde{P}}_2+\sum\limits_{i=1}^3 a_i^2 \check{\tilde{P}}_2+Q_2+\frac{\bar{B}_1 \check{\tilde{P}}_2 \tilde{P}}{R_1+\left(\sum\limits_{i=1}^3 b_i^2\right)\tilde{P}}=0, \quad
\check{\tilde{P}}_2(T)=G_2,
\end{equation}
where $\tilde{P}(\cdot)$ is the solution to (\ref{added Riccati equations of followers}). Since ODEs (\ref{auxiliary ODEs 1 of followers}) and (\ref{auxiliary ODEs 2 of followers}) are both linear, they have unique solutions $\check{\tilde{P}}_1(\cdot)$ and $\check{\tilde{P}}_2(\cdot)$, respectively. On the other hand, it is easy to see that $\tilde{P}_1(\cdot)$ and $\tilde{P}_2(\cdot)$ in (\ref{coupled Riccati equations of followers}) also solutions to (\ref{auxiliary ODEs 1 of followers}) and (\ref{auxiliary ODEs 2 of followers}), respectively. It follows from the uniqueness of solutions to (\ref{auxiliary ODEs 1 of followers}) and (\ref{auxiliary ODEs 2 of followers}) that
$$
\check{\tilde{P}}_1(\cdot)\equiv\tilde{P}_1(\cdot), \quad \check{\tilde{P}}_2(\cdot)\equiv\tilde{P}_2(\cdot) .
$$
Then it means that system (\ref{coupled Riccati equations of followers}) has the unique solutions $\tilde{P}_1(\cdot)$ and $\tilde{P}_2(\cdot)$.
\end{proof}

{\it Proof of Theorem \ref{thm3}.}\begin{proof}
{\it Necessity.} We introduce the leaders' Hamiltonian functions as, for $j=3,4$,
\begin{equation*}\label{Hamilton of leaders}
\begin{aligned}
&H_j\big(t, x^{u_3, u_4}, u_3, u_4, \hat{\Phi}, \hat{\zeta}_1, \hat{\zeta}_2 ; y_j, z_j, q_{j 1}, q_{j 2}, q_{j 3}\big) \\
& := z_j\big(a_0 x^{u_3, u_4}+M_{01} \hat{x}^{u_3, u_4} +N_{00}^T \hat{\Phi} +N_{01}^T \hat{\zeta}_1 +N_{03}^T \hat{\zeta}_3 +d_0 u_3+M_{04} \hat{u}_3\big) \\
&\quad +e_0 u_4+M_{05} \hat{u}_4+\sum\limits_{i=1}^3 q_{ji}\big( a_i x^{u_3 u_4}+M_{i 1} \hat{x}^{u_3 u_4}+N_{i 0}^T \hat{\Phi} +N_{i 1}^T \hat{\zeta_1} +N_{i 3}^T \hat{\zeta}_3  \\
&\quad +d_i u_3+M_{i 4} \hat{u}_3+e_i u_4+M_{i 5} \hat{u}_4\big) +\frac{1}{2} \Big[Q_j (x^{u_3, u_4})^2+R_j u_j^2 \Big] \\
&\quad +\big\langle y_j,A_0 \hat{\Phi}+ A_1 \hat{\zeta}_1+A_3 \hat{\zeta}_3+f_1 \hat{u}_3  +f_2 \hat{u}_4\big\rangle.
\end{aligned}
\end{equation*}
It follows from Proposition 2.3 in \cite{Shi-Wang-Xiong-2016} that, if $\left(u_3^*(\cdot), u_4^*(\cdot)\right)$ is a Nash equilibrium point of the leaders, then (\ref{SMP of leaders}) holds.

{\it Sufficiency.} Next, we prove that $\left(u_3^*(\cdot), u_4^*(\cdot)\right)$ in (\ref{SMP of leaders}) is indeed a Nash equilibrium point of the leaders. For any $u_3(\cdot) \in \mathcal{U}_3$, letting $x_3(\cdot):=x^{u_3, u_4^*}(\cdot)-x^*(\cdot)$, $\hat{\Phi}_3(\cdot):=\hat{\Phi}(\cdot)-\hat{\Phi}^*(\cdot)$, $\hat{\tilde{\zeta}}_1(\cdot):=\hat{\zeta}_1(\cdot)-\hat{\zeta}_1^*(\cdot)$, $\hat{\tilde{\zeta}}_3(\cdot):=\hat{\zeta}_3(\cdot)-\hat{\zeta}_3^*(\cdot)$, we have
$$
\hat{J}_3\left(u_3(\cdot), u_4^*(\cdot) \right)-\hat{J}_3\left(u_3^*(\cdot), u_4^*(\cdot) \right)
=\frac{1}{2} \mathbb{E}\left[\int_0^T \Big(Q_3 x_3^2+R_3\left(u_3-u_3^*\right)^2\Big) d t+G_3 x_3^2(T)\right]+\Theta_l,
$$
where
$$
\Theta_l:= \mathbb{E}\left[\int_0^T \Big(Q_3 x_3 x^*+R_3\left(u_3-u_3^*\right) u_3^*\Big) d t+G_3 x_3(T) x^*(T)\right].
$$
Notice that
$$
\mathbb{E}\big[G_3 x_3(T) x^*(T)\big]= \mathbb{E}\big[x_3(T) z_3(T)\big],\quad x_3(0)=0=\mathbb{E}\big\langle y_3(0),\hat{\Phi}(0)\big\rangle = \mathbb{E}\big\langle y_3(T),\hat{\Phi}(T)\big\rangle.
$$
Applying It\^{o}'s formula to $x_3(\cdot) z_3(\cdot))$ and $\big\langle y_3(\cdot), \hat{\Phi}(\cdot)\big\rangle $ respectively, we get
\begin{equation*}
\begin{aligned}
& \mathbb{E} \big[ x_3(T) z_3(T)-x_3(0) z_3(0)\big]\\
& =\mathbb{E} \int_0^T\bigg\{N_{00}^\top \hat{\Phi} z_3+ N_{01}^\top \hat{\tilde{\zeta}}_1 z_3+ N_{03}^\top \hat{\tilde{\zeta}}_3 z_3
 +d_0\left(u_3-u_3^*\right)z_3+M_{04}\left(\hat{u}_3-\hat{u}_3^* \right) z_3- Q_3 x_3 x^* \\
&\qquad +\sum_{i=1}^3\Big[N_{i0}^\top \hat{\Phi} q_{3i} +N_{i1}^\top \hat{\tilde{\zeta}}_1 q_{3i}\ + N_{i3}^\top \hat{\tilde{\zeta}}_3 q_{3i}+d_i\left(u_3-u_3^*\right)q_{3 i}
 +M_{i 4}\left(\hat{u}_3-\hat{u}_3^*\right) q_{3 i}\Big]\bigg\} d t,
\end{aligned}
\end{equation*}
and
\begin{equation*}
\begin{aligned}
0& = - \mathbb{E} \big[ \big\langle y_3(T), \hat{\Phi}(T)\big\rangle+\big\langle y_3(0), \hat{\Phi}(0)\big\rangle\big] =\mathbb{E} \int_0^{T}\bigg[\big\langle f_1(\hat{u}_3-\hat{u}_3^*),y_3\big\rangle\\
&\quad -\bigg\langle N_{00} p_3 + \sum\limits_{i=1}^3 N_{i0} q_{3i}, \hat{\Phi}\bigg\rangle
 -\bigg\langle N_{01}z_3+ \sum\limits_{i=1}^3 N_{i1} q_{3i}, \hat{\zeta}_1\bigg\rangle -\big\langle N_{03}z_3,\hat{\zeta}_3\big\rangle \bigg] dt.
\end{aligned}
\end{equation*}
By (\ref{SMP of leaders}), we obtain
\begin{equation*}
\begin{aligned}
\Theta_l& =\mathbb{E} \int_0^T\bigg[\left(u_3-u_3^*\right) \bigg(R_3^*+ d_0 z_3+\sum_{i=1}^3 d_i q_{3i}\bigg) \\
&\qquad\qquad +\left(\hat{u}_3-\hat{u}_3^*\right) \bigg(M_{04}^* z_3+ \sum_{i=1}^3 M_{i4} q_{3i} +\left\langle f_1, y_3\right\rangle \bigg) \bigg] dt\\
& =\mathbb{E} \int_0^T \left\{\mathbb{E}\bigg[\left(u_3-u_3^*\right)\bigg(R_3 u_3^*+d_0 z_3+\sum_{i=1}^3 d_i q_{3 i}\bigg)\right.\\
&\qquad\qquad \left.+\left(\hat{u}_3-\hat{u}_3^*\right)\bigg(M_{04} z_3+\sum_{i=1}^3 M_{i 4} q_{3 i}+\left\langle f_1, y_3 \right\rangle \bigg) \bigg| \mathcal{G}_t^2 \bigg]\right\} d t\\
\end{aligned}
\end{equation*}
\begin{equation*}
\begin{aligned}
& =\mathbb{E} \int_0^T \mathbb{E}\bigg[\left(u_3-u_3^*\right)\bigg(-M_{04} \check{\hat{z}}_3-\sum_{j=1}^3 M_{j 4} \check{\hat{q}}_{3 j}- f_1^\top \check{\hat{y}}_3 \bigg) \\
&\qquad\qquad +\left(\hat{u}_3-\hat{u}_3^*\right)\bigg(M_{04} \hat{z}_3+\sum_{i=1}^3 M_{i 4} \check{q}_{3 i}+\left\langle f_1, \check{y}_3 \right\rangle \bigg)  \bigg] d t = 0.
\end{aligned}
\end{equation*}
Here we used the property of the conditional expectation, for several times.

Since Assumption \ref{A4} holds, we get
$$
\hat{J}_3\left( u_3(\cdot), u_4^*(\cdot)\right)-\hat{J}_3\left( u_3^*(\cdot), u_4^*(\cdot)\right) \geq 0 .
$$
Similarly, for any $u_4(\cdot) \in \mathcal{U}_4$, we also have
$$
\hat{J}_4\left( u_3^*(\cdot), u_4(\cdot)\right)-\hat{J}_4\left( u_3^*(\cdot), u_4^*(\cdot)\right) \geq 0 .
$$
Therefore, by (\ref{Nash game of leaders}), $\left(u_3^*(\cdot), u_4^*(\cdot)\right)$ is a Nash equilibrium point for the leaders' problem.
\end{proof}

{\it Derivation of (\ref{Z--}).}\begin{proof}
First, applying first $\mathbb{E}\left[\mathbb{E} \left[ \cdot \mid \mathcal{G}_t^1\right] \mid \mathcal{G}_t^2 \right]$ to (\ref{coupled Z1}), (\ref{coupled Z2}) and (\ref{coupled Z3}), we can get
\begin{equation}\label{hat and check Z}
\check{\hat{Z}}_i=K_{i 0} \check{\hat{X}}+K_{i 1} \check{\hat{Z}}_1+K_{i 2} \check{\hat{Z}}_2+K_{i 3} \check{\hat{Z}}_3, \quad i=1,2,3,
\end{equation}
where
$$
\begin{aligned}
K_{10}&:=\left(P_1+P_2\right)\left(\mathcal{A}_1+\mathcal{M}_{11}+\mathcal{H}_1\right)+\left(P_1+P_2\right)\left(\mathcal{N}_{10}+\mathcal{B}_{10}+\mathcal{C}_0\right)\left(P_1+P_2+P_3+P_4\right), \\
K_{20}&:=\left(P_1+P_3\right)\left(\mathcal{A}_2+\mathcal{M}_{21}+\mathcal{H}_2\right)+\left(P_1+P_3\right)\left(\mathcal{N}_{20}+\mathcal{B}_{20}+\mathcal{C}_{20}\right)\left(P_1+P_2+P_3+P_4\right), \\
K_{30}&:=\left(P_1+P_2+P_3+P_4\right)\left(\mathcal{A}_3+\mathcal{M}_{31}+\mathcal{H}_3\right)+\left(P_1+P_2+P_3+P_4\right)^2 \left(\mathcal{N}_{30}+\mathcal{B}_{30}+\mathcal{C}_{30}\right), \\
K_{1i}&:=\left(P_1+P_2\right)\left(\mathcal{N}_{1i}+\mathcal{B}_{1i}+\mathcal{C}_{1i}\right), \quad K_{2i}:=\left(P_1+P_3\right)\left(\mathcal{N}_{2i}+\mathcal{B}_{2i}+\mathcal{C}_{2i}\right), \\
K_{3i}&:=\left(P_1+P_2+P_3+P_4\right)\left(\mathcal{N}_{3i}+\mathcal{B}_{3i}+\mathcal{C}_{3i}\right), \quad i=1,2,3.
\end{aligned}
$$
Therefore, (\ref{hat and check Z}) can be rewritten as
\begin{equation}\label{Equations of hat and check Z}
\left(\begin{array}{ccc}
I_4-K_{11} & -K_{12} & -K_{13} \\
-K_{21} & I_4-K_{22} & -K_{23} \\
-K_{31} & -K_{32} & I_4-K_{33}
\end{array}\right)\left(\begin{array}{c}
\check{\hat{Z}}_1 \\
\check{\hat{Z}}_2 \\
\check{\hat{Z}}_3
\end{array}\right)=\left(\begin{array}{l}
K_{10} \\
K_{20} \\
K_{30}
\end{array}\right) \check{\hat{X}}.
\end{equation}
We need the following assumption.
\begin{assumption}\label{A5}
The coefficient matrix of (\ref{Equations of hat and check Z}) is invertible, for any $t \in \left[0,T\right]$.
\end{assumption}
Under this, we can derive that
\begin{equation}\label{hat and check Z-}
\check{\hat{Z}}_i=\mathbb{N}_1^{-1}\left(\textbf{K}_{1 i} K_{10}+\textbf{K}_{2 i} K_{20}+\textbf{K}_{3 i} K_{30}\right) \check{\hat{X}}\equiv N_i \check{\hat{X}}, \quad i=1,2,3,
\end{equation}
where $\mathbb{N}_1$ is the determinant of the coefficient matrix of (\ref{Equations of hat and check Z}), and $\textbf{K}_{j i}$ is the the $(i,j)$ block matrix in adjoint matrix of the above coefficient matrix, for $j, i = 1, 2, 3$.

Next, applying $\mathbb{E} \left[ \cdot \mid \mathcal{G}_t^2\right] $ to (\ref{coupled Z1}), (\ref{coupled Z2}) and (\ref{coupled Z3}), we can obtain
\begin{equation}\label{check Z}
\check{Z}_i=\overline{K}_{i 0} \check{X}+ \tilde{K}_{i 0} \check{\hat{X}}+\overline{K}_{i 1} \check{Z}_1+\tilde{K}_{i 1} \check{\hat{Z}}_1+\overline{K}_{i 2} \check{Z}_2+\tilde{K}_{i 2}\check{\hat{Z}}_2 +\overline{K}_{i 3} \check{Z}_3+\tilde{K}_{i 3}\check{\hat{Z}}_3, \quad i=1,2,3,
\end{equation}
where
$$
\begin{aligned}
\overline{K}_{10} &:=P_1 \mathcal{A}_1+P_1\left(\mathcal{N}_{10}+\mathcal{B}_{10}\right)\left(P_1+P_3\right), \\
\tilde{K}_{10} &:= P_1\left(\mathcal{M}_{11}+\mathcal{H}_1\right)+P_1\left(\mathcal{N}_{10}+\mathcal{B}_{10}\right)\left(P_2+P_4\right)+P_2\left(\mathcal{A}_1+\mathcal{M}_{11}+\mathcal{H}_1\right) \\
&\quad +\big[P_1 \mathcal{C}_{10}+P_2\left(\mathcal{N}_{10}+\mathcal{B}_{10}+\mathcal{C}_{10}\right)\big]\left(P_1+P_2+P_3+P_4\right), \\
\overline{K}_{20} &:= \left(P_1+P_3\right) \mathcal{A}_2+\left(P_1+P_3\right)\left(\mathcal{N}_{20}+\mathcal{B}_{20}\right)\left(P_1+P_3\right), \\
\tilde{K}_{20} &:= \left(P_1+P_3\right)\left(\mathcal{M}_{21}+\mathcal{H}_2\right)+\left(P_1+P_3\right)\left(\mathcal{N}_{20}+\mathcal{B}_{20}\right)\left(P_2+P_4\right) \\
&\quad +\left(P_1+P_3\right) \mathcal{C}_{20}\left(P_1+P_2+P_3+P_4\right), \\
\overline{K}_{30} &:=\left(P_1+P_3\right) \mathcal{A}_3+\left(P_1+P_3\right)\left(\mathcal{N}_{30}+\mathcal{B}_{30}\right)\left(P_1+P_3\right), \\
 \tilde{K}_{30} &:=\left(P_1+P_3\right)\left(\mathcal{M}_{31}+\mathcal{H}_3\right)+\left(P_1+P_3\right)\left(\mathcal{N}_{30}+\mathcal{B}_{30}\right)\left(P_2+P_4\right)
 +\left(P_2+P_4\right)\left(\mathcal{A}_3+\mathcal{M}_{31}\right. \\
&\quad \left. +\mathcal{H}_3\right)+\big[\left(P_1+P_3\right) \mathcal{C}_{30}+\left(P_2+P_4\right)\left(\mathcal{B}_{30}+\mathcal{C}_{30}+\mathcal{N}_{30}\right)\big]\left(P_1+P_2+P_3+P_4\right), \\
\overline{K}_{1i} &:= P_1\left(\mathcal{N}_{1i}+\mathcal{B}_{1i}\right), \ \tilde{K}_{1i}:=P_1 \mathcal{C}_{1i}+P_2\left(\mathcal{N}_{1i}+\mathcal{B}_{1i}+\mathcal{C}_{1i}\right),\
\overline{K}_{2i} := \left(P_1+P_3\right)\left(\mathcal{N}_{2i}+\mathcal{B}_{2i}\right), \\
\tilde{K}_{2i} &:=\left(P_1+P_3 \right)\mathcal{C}_{2i}, \ \overline{K}_{3i}:= \left(P_1+P_3\right)\left(\mathcal{N}_{3i}+\mathcal{B}_{3i}\right), \\
\tilde{K}_{3i} &:=\left(P_1+P_3\right) \mathcal{C}_{3i}+\left(P_2+P_4\right)\left(\mathcal{N}_{3i}+\mathcal{B}_{3i}+\mathcal{C}_{3i}\right), \quad i=1,2,3.
\end{aligned}
$$
Substituting (\ref{hat and check Z-}) into (\ref{check Z}), we obtain
\begin{equation}\label{check Z-}
\begin{aligned}
\check{Z}_i &=\overline{K}_{i 0} \check{X}+\big(\tilde{K}_{i 0}+\tilde{K}_{i 1} N_1+\tilde{K}_{i 2} N_2+\tilde{K}_{i 3} N_3\big) \check{\hat{X}}+\overline{K}_{i 1} \check{Z}_1
+\overline{K}_{i 2} \check{Z}_2+\overline{K}_{i 3} \check{Z} \\
& \equiv \overline{K}_{i 0} \check{X}+\overline{N}_i \check{\hat{X}}+\overline{K}_{i 1} \check{Z}_1+\overline{K}_{i 2} \check{Z}_2+\overline{K}_{i 3} \check{Z}_3, \quad i=1,2,3.
\end{aligned}
\end{equation}
It can be rewritten as
\begin{equation}\label{Equations of check Z}
\left(\begin{array}{ccc}
I_4-\overline{K}_{11} & -\overline{K}_{12} & -\overline{K}_{13} \\
-\overline{K}_{21} & I_4-\overline{K}_{22} & -\overline{K}_{23} \\
-\overline{K}_{31} & -\overline{K}_{32} & I_4-\overline{K}_{33}
\end{array}\right)\left(\begin{array}{c}
\check{Z}_1 \\
\check{Z}_2 \\
\check{Z}_3
\end{array}\right)=\left(\begin{array}{l}
K_{10}\check{X}_t+\overline{N}_1 \check{\hat{X}}\\
K_{20}\check{X}_t+\overline{N}_2 \check{\hat{X}}\\
K_{30}\check{X}_t+\overline{N}_3 \check{\hat{X}}
\end{array}\right).
\end{equation}
Here, we need the following assumption.
\begin{assumption}\label{A6}
The coefficient matrix of (\ref{Equations of check Z}) is invertible, for any $t \in \left[0,T\right]$.
\end{assumption}
Then, we can derive that
\begin{equation}\label{check Z--}
\begin{aligned}
\check{Z}_i &=\mathbb{N}_2^{-1}\left[\left(\check{\textbf{K}}_{1 i} \overline{K}_{10}+\check{\textbf{K}}_{2 i} \overline{K}_{20}+\check{\textbf{K}}_{3i} \overline{K}_{30}\right) \check{X}+\left(\check{\textbf{K}}_{1 i} \overline{N}_{1}+\check{\textbf{K}}_{2 i} \overline{N}_{2}+\check{\textbf{K}}_{3i} \overline{N}_{3}\right)\check{\hat{X}} \right] \\
&\equiv \check{N}_i \hat{X}+\tilde{N}_i \check{\hat{X}}, \quad i=1,2,3,
\end{aligned}
\end{equation}
where $\mathbb{N}_2$ is the determinant of the coefficient matrix of (\ref{Equations of check Z}), and $\check{\textbf{K}}_{j i}$ is the the $(i,j)$ block matrix in adjoint matrix of the above coefficient matrix, for $j, i = 1, 2, 3$.

In the next step, applying $\mathbb{E} \left[ \cdot \mid \mathcal{G}_t^1\right] $ to (\ref{coupled Z1}), (\ref{coupled Z2}) and (\ref{coupled Z3}), we can obtain
\begin{equation}\label{hat Z}
\hat{Z}_i=\hat{K}_{i 0} \hat{X}+\check{\hat{K}}_{i 0} \check{\hat{X}}+\hat{K}_{i 1} \hat{Z}_1+\check{\hat{K}}_{i 1} \check{\hat{Z}}_1+\hat{K}_{i 2} \hat{Z}_2+\check{\hat{X}}_{i 2} \check{\hat{X}}_2+\hat{K}_{i 3} \hat{Z}_3+\check{\hat{K}}_{i 3} \check{\hat{Z}}_3, \quad i=1,2,3,
\end{equation}
where
$$
\begin{aligned}
\hat{K}_{10} &:= \left(P_1+P_2\right)\left(\mathcal{A}_1+\mathcal{M}_{11}\right)+\left(P_1+P_2\right) \mathcal{N}_{10}\left(P_1+P_2\right), \\
\check{\hat{K}}_{10} &:= \left(P_1+P_2\right) \mathcal{N}_{10}\left(P_3+P_4\right)+\left(P_1+P_2\right) \mathcal{H}_1+\left(P_1+P_2\right)\left(\mathcal{B}_{10}+\mathcal{C}_{10}\right)\left(P_1+P_2+P_3+P_4\right), \\
\hat{K}_{20} &:= P_1\left(\mathcal{A}_2+\mathcal{M}_{21}\right)+P_1 \mathcal{N}_{20}\left(P_1+P_2\right), \\
\check{\hat{K}}_{20} &:= P_1 \mathcal{H}_2+P_3\left(\mathcal{A}_2+\mathcal{M}_{21}+\mathcal{H}_2\right)+P_1 \mathcal{N}_{20}\left(P_3+P_4\right) \\
&\quad +\left[P_1\left(\mathcal{B}_{20}+\mathcal{C}_{20}\right)+P_3\left(\mathcal{N}_{20}+\mathcal{B}_{20}+\mathcal{C}_{20}\right)\right]\left(P_1+P_2+P_3+P_4\right) \\
\hat{K}_{30} &:= \left(P_1+P_2\right)\left(\mathcal{A}_3+\mathcal{M}_{31}\right)+\left(P_1+P_2\right) \mathcal{N}_{30}\left(P_1+P_2\right), \\
\check{\hat{K}}_{30} &:= \left(P_1+P_2\right) \mathcal{H}_3+\left(P_3+P_4\right)\left(\mathcal{A}_3+\mathcal{M}_{31}+\mathcal{H}_3\right)+\left(P_1+P_2\right) \mathcal{N}_{30}\left(P_3+P_4\right) \\
&\quad +\left[\left(P_1+P_2\right)\left(\mathcal{B}_{30}+\mathcal{C}_{30}\right)+\left(P_3+P_4\right)\left(\mathcal{N}_{30}+\mathcal{B}_{30}+\mathcal{C}_{30}\right)\right]\left(P_1+P_2+P_3+P_4\right), \\
\hat{K}_{1i} &:= \left(P_1+P_2\right)\mathcal{N}_{1i}, \quad \check{\hat{K}}_{1i}:=\left(P_1+P_2\right) \left(\mathcal{B}_{1i}+\mathcal{C}_{1i}\right),\quad \hat{K}_{2i} := P_1 \mathcal{N}_{2i}, \\
\check{\hat{K}}_{2i} &:= P_1 \left(\mathcal{B}_{2i}+\mathcal{C}_{2i}\right)+P_3\left(\mathcal{N}_{2i}+\mathcal{B}_{2i}+\mathcal{C}_{2i}\right),\quad \hat{K}_{3i} := \left(P_1+P_2\right) \mathcal{N}_{3i}, \\
\check{\hat{K}}_{3i} &:=\left(P_1+P_2\right) \left(\mathcal{B}_{3i}+\mathcal{C}_{3i}\right)+\left(P_3+P_4\right)\left(\mathcal{N}_{3i}+\mathcal{B}_{3i}+\mathcal{C}_{3i}\right), \quad i=1,2,3.
\end{aligned}
$$
Substituting (\ref{hat and check Z}) into (\ref{hat Z}), we obtain
\begin{equation}\label{hat Z-}
\begin{aligned}
\hat{Z}_i & =\hat{K}_{i 0} \hat{X}+\big(\check{\hat{K}}_{i 0}+\check{\hat{K}}_{i 1} N_1+\check{\hat{K}}_{i 2} N_2+\check{\hat{K}}_{i 3} N_3\big) \check{\hat{X}}+\hat{K}_{i 1} \hat{Z}_1+\hat{K}_{i 2} \hat{Z}_2+\hat{K}_{i 3} \hat{Z}_3 \\
& \equiv \hat{K}_{i 0} \hat{X}+\overline{\overline{N}}_i \check{\hat{X}}+\hat{K}_{i 1} \hat{Z}_1+\hat{K}_{i 2} \hat{Z}_2+\hat{K}_{i 3} \hat{Z}_3,\quad i=1,2,3.
\end{aligned}
\end{equation}
It can also be rewritten as
\begin{equation}\label{Equations of hat Z}
\left(\begin{array}{ccc}
I_4-\hat{K}_{11} & -\hat{K}_{12} & -\hat{K}_{13} \\
-\hat{K}_{21} & I_4-\hat{K}_{22} & -\hat{K}_{23} \\
-\hat{K}_{31} & -\hat{K}_{32} & I_4-\hat{K}_{33}
\end{array}\right)\left(\begin{array}{c}
\hat{Z}_t^1 \\
\hat{Z}_t^2 \\
\hat{Z}_t^3
\end{array}\right)=\left(\begin{array}{l}
\hat{K}_{10}\hat{X}_t+\overline{\overline{N}}_1 \check{\hat{X}}_t\\
\hat{K}_{20}\hat{X}_t+\overline{\overline{N}}_2 \check{\hat{X}}_t\\
\hat{K}_{30}\hat{X}_t+\overline{\overline{N}}_3 \check{\hat{X}}_t
\end{array}\right)
\end{equation}
Additionally, we need the following assumption.
\begin{assumption}\label{A7}
The coefficient matrix of (\ref{Equations of hat Z}) is invertible, for any $t \in \left[0,T\right]$.
\end{assumption}
Then, we can get that
\begin{equation}\label{hat Z--}
\begin{aligned}
\hat{Z}_i &={\mathbb{N}_3^{-1}}\big[\big(\hat{\textbf{K}}_{1 i} \hat{{K}}_{10}+\hat{\textbf{K}}_{2 i} \hat{K}_{20}+\hat{\textbf{K}}_{3i} \hat{K}_{30}\big) \hat{X}+\big(\hat{\textbf{K}}_{1 i} \overline{\overline{N}}_{1}+\hat{\textbf{K}}_{2 i} \overline{\overline{N}}_{2}+\hat{\textbf{K}} \overline{\overline{N}}_{3}\big)\check{\hat{X}} \big] \\
&\equiv \hat{N}_i \hat{X}+\check{\hat{N}}_i \check{\hat{X}}, \quad i=1,2,3,
\end{aligned}
\end{equation}
where $\mathbb{N}_3$ is the determinant of the coefficient matrix of (\ref{Equations of hat Z}), and $\hat{\textbf{K}}_{j i}$ is the the $(i,j)$ block matrix in adjoint matrix of the above coefficient matrix, for $j,i =1,2,3$.

In final step, substituting (\ref{hat and check Z-}), (\ref{check Z--}) and (\ref{hat Z--}) into (\ref{coupled Z1}), (\ref{coupled Z2}) and (\ref{coupled Z3}), we obtain
\begin{equation}\label{Z}
\begin{aligned}
Z_i & =\mathcal{K}_{i 0} X+\hat{\mathcal{K}}_{i 0} \hat{X}+\check{\mathcal{K}}_{i 0} \check{X}+\check{\hat{\mathcal{K}}}_{i 0} \check{\hat{X}} +P_1 \mathcal{N}_{i1} Z_1 +P_1 \mathcal{N}_{i2} Z_2 +P_1\mathcal{N}_{i3} Z_3, \quad i=1,2,3,
\end{aligned}
\end{equation}
where
$$
\begin{aligned}
\mathcal{K}_{1i} &:=P_1 \mathcal{A}_i+P_1^2 \mathcal{N}_{i0}, \quad i=1,2,3, \\
\hat{\mathcal{K}}_{10} &:=P_2 \mathcal{A}_1+\left(P_1+P_2\right) \mathcal{M}_{11}+P_1 \mathcal{N}_{10} P_2+P_2 \mathcal{N}_{10}\left(P_1+P_2\right)+P_2 \sum_{i=1}^3 \mathcal{N}_{1i} \hat{N}_i, \\
\check{\mathcal{K}}_{10} &:=P_1 \mathcal{N}_{10} P_3+P_1 \mathcal{B}_{10}\left(P_1+P_3\right)+P_1 \sum_{i=1}^3 \mathcal{B}_{1i} \tilde{N}_i, \\
\check{\hat{\mathcal{K}}}_{10} &:=\left(P_1+P_2\right) \mathcal{H}_1+P_1 \mathcal{N}_{10} P_4+P_2 \mathcal{N}_{10}\left(P_3+P_4\right)+P_1 \mathcal{B}_{10}\left(P_2+P_4\right)\\
&\quad +\big[P_1 \mathcal{C}_{10}+P_2\left(\mathcal{B}_{10}+\mathcal{C}_{10}\right)\big]\left(P_1+P_2+P_3+P_4\right) \\
&\quad +\sum_{i=1}^3\left[P_2 \mathcal{N}_{1i} \check{\hat{N}}_i+P_1 \mathcal{B}_{1i} \tilde{N}_i+\left(P_1 \mathcal{C}_{1 i}+P_2 \mathcal{B}_{1i}+P_2 \mathcal{C}_{1 i}\right) N_i\right] \\
\hat{\mathcal{K}}_{20} &:=P_1 \mathcal{M}_{21}+P_1 \mathcal{N}_{20} P_2, \quad \check{K}_{20} :=P_1 \mathcal{N}_{20} P_3+P_1 \mathcal{B}_{20}\left(P_1+P_3\right)+P_3 \mathcal{A}_2\\
&\quad +P_3\left(\mathcal{N}_{20}+\mathcal{B}_{20}\right)\left(P_1+P_3\right)+\sum_{i=1}^3\left(P_1 \mathcal{B}_{2 i}+P_3 \mathcal{N}_{2 i}+P_3 \mathcal{B}_{2 i}\right) \tilde{N}_i, \\
\check{\hat{\mathcal{K}}}_{20} &:=P_1 \mathcal{H}_2+P_1  \mathcal{N}_{20} P_4+P_1 \mathcal{B}_{20}\left(P_2+P_4\right)+P_1 \mathcal{C}_{20}\left(P_1+P_2+P_3+P_4\right) \\
&\quad +P_3\left(\mathcal{M}_{21}+\mathcal{H}_2\right)+P_3\left(\mathcal{N}_{20}+\mathcal{B}_{20}\right)\left(P_2+P_4\right)+P_3 \mathcal{C}_{20}\left(P_1+P_2+P_3+P_4\right) \\
&\quad +\sum_{i=1}^3\left[\left(P_1 B_{2 i}+P_3 N_{2 i}+P_3 B_{2 i}\right) \tilde{N}_i+\left(P_1 C_{2 i}+P_3 C_{2 i}\right) N_i\right], \\
\hat{\mathcal{K}}_{30} &:=P_1 \mathcal{M}_{31}+P_1 \mathcal{N}_{30} P_2+P_2\left(\mathcal{A}_3+\mathcal{M}_{31}\right)+P_3 \mathcal{N}_{30}\left(P_1+P_2\right)+\sum_{i=1}^3 P_2 \mathcal{N}_{3 i} \hat{N}_i, \\
\check{\mathcal{K}}_{30} &:=P_1 \mathcal{N}_{30} P_3+P_1 \mathcal{B}_{30}\left(P_1+P_3\right)+P_3 \mathcal{A}_3+P_3 \left( \mathcal{N}_{30}+\mathcal{B}_{30}\right)\left(P_1+P_3\right) \\
&\quad +\sum_{i=1}^3\left(P_1 \mathcal{B}_{3 i}+P_3 \mathcal{N}_{3 i}+P_3 \mathcal{B}_{3 i}\right) \overline{N}_i, \\
\check{\hat{\mathcal{K}}}_{30} &:=\left(P_1+P_2\right) \mathcal{H}_3+P_3\left(\mathcal{M}_{31}+\mathcal{H}_3\right)+P_4\left(\mathcal{A}_3+\mathcal{M}_{31}+\mathcal{H}_3\right)+P_1 \mathcal{NN}_{30} P_4+P_1 \mathcal{B}_{30}\left(P_2+P_4\right) \\
&\quad +P_2 \mathcal{N}_{30}\left(P_3+P_4\right)+P_3\left(\mathcal{N}_{30}+\mathcal{B}_{30}\right)\left(P_2+P_4\right)+\big[\left(P_1+P_3\right) \mathcal{C}_{30}+P_2\left(\mathcal{B}_{30}+\mathcal{C}_{30}\right) \\
&\quad +P_4\left(\mathcal{N}_{30}+\mathcal{B}_{30}+\mathcal{C}_{30}\right)\big]\left(R_1+P_2+P_3+P_4\right)+\sum_{i=1}^3\left\{\left(P_1 \mathcal{B}_{3 i}+P_3 \mathcal{N}_{3 i}+P_3 \mathcal{B}_{3 i}\right) \tilde{N_i}\right. \\
&\quad \left. +P_2 \mathcal{N}_{3 i} \check{\hat{N}}_i+\big[\left(P_1+P_3\right) \mathcal{C}_{3 i}+P_2\left(\mathcal{B}_{3 i}+\mathcal{C}_{3 i}\right)+P_4\left(\mathcal{N}_{3 i}+\mathcal{B}_{3 i}+\mathcal{C}_{3 i}\right)\big] N_i\right\}.
\end{aligned}
$$
To solve (\ref{Z}), we rewrite it as
\begin{equation}\label{Equations of Z}
\left(\begin{array}{ccc}
I_4-P_1 \mathcal{N}_{11} & -P_1 \mathcal{N}_{12} & -P_1 \mathcal{N}_{13} \\
-P_{10} \mathcal{N}_{21} & I_4-P_1 \mathcal{N}_{22} & -P_1 \mathcal{N}_{23} \\
-P_1 \mathcal{N}_{31} & -P_1 \mathcal{N}_{32} & I_4- P_1 \mathcal{N}_{33}
\end{array}\right)\left(\begin{array}{l}
Z_1 \\
Z_2 \\
Z_3
\end{array}\right)=\left(\begin{array}{l}
\mathcal{K}_{10} X+\hat{\mathcal{K}}_{10} \hat{X}+\check{\mathcal{K}}_{10} \check{X}+\check{\hat{\mathcal{K}}}_{10} \check{\hat{X}} \\
\mathcal{K}_{20} X+\hat{\mathcal{K}}_{20} \hat{X}+\check{\mathcal{K}}_{20} \check{X}+\check{\hat{\mathcal{K}}}_{20} \check{\hat{X}} \\\
\mathcal{K}_{30} X+\hat{\mathcal{K}}_{30} \hat{X}+\check{\mathcal{K}}_{30} \check{X}+\check{\hat{\mathcal{K}}}_{30} \check{\hat{X}} \\
\end{array}\right).
\end{equation}
We also need the following additional assumption.
\begin{assumption}\label{A8}
The coefficient matrix of (\ref{Equations of Z}) is invertible, for any $t \in \left[0,T\right]$.
\end{assumption}
Then, we can obtain that the following state-estimate feedback form of $Z_i$:
\begin{equation}\label{Z-}
\begin{aligned}
Z_i &={\mathbb{N}_4^{-1}}\left[\big(\check{\hat{\textbf{K}}}_{1i} \mathcal{K}_{10}+\check{\hat{\textbf{K}}}_{2 i} \mathcal{K}_{20}+\check{\hat{\textbf{K}}}_{3i} \mathcal{K}_{30}\big) X
 +\big(\check{\hat{\textbf{K}}}_{1i} \hat{\mathcal{K}}_{10}+\check{\hat{\textbf{K}}}_{2 i} \hat{\mathcal{K}}_{20}+\check{\hat{\textbf{K}}}_{3i} \hat{\mathcal{K}}_{30}\big) \hat{X} \right.\\
&\quad \left.+\big(\check{\hat{\textbf{K}}}_{1i} \check{\mathcal{K}}_{10}+\check{\hat{\textbf{K}}}_{2 i} \check{\mathcal{K}}_{20}+\check{\hat{\textbf{K}}}_{3i} \check{\mathcal{K}}_{30}\big) \check{X}
 +\big(\check{\hat{\textbf{K}}}_{1i} \check{\hat{\mathcal{K}}}_{10}+\check{\hat{\textbf{K}}}_{2 i} \check{\hat{\mathcal{K}}}_{20}+\check{\hat{\textbf{K}}}_{3i} \check{\hat{\mathcal{K}}}_{30}\big) \check{\hat{X}} \right] \\
&\equiv \mathcal{N}_i X+\hat{\mathcal{N}}_i \hat{X} +\check{\mathcal{N}}_i \check{X}+\check{\hat{\mathcal{N}}}_i \check{\hat{X}}, \quad i=1,2,3,
\end{aligned}
\end{equation}
where $\mathbb{N}_4$ is the determinant of the coefficient matrix of (\ref{Equations of Z}), and $\check{\hat{\textbf{K}}}_{j i}$ is the the $(i,j)$ block matrix in adjoint matrix of the above coefficient matrix, for $j, i = 1, 2, 3$.
\end{proof}

\end{document}